\documentclass{article}
\usepackage[english]{babel}
\usepackage[top=2cm,bottom=2cm,left=3cm,right=3cm,marginparwidth=1.75cm]{geometry}
\usepackage{amsmath}
\usepackage[utf8]{inputenc}   
\usepackage[T1]{fontenc}      
\usepackage[english]{babel} 
\usepackage{amssymb}
\usepackage{mathrsfs}
\usepackage{graphicx}
\usepackage[utf8]{inputenc}
\usepackage{newunicodechar}
\newunicodechar{́}{\'}
\usepackage[colorlinks=true, allcolors=blue]{hyperref}
\usepackage[T1]{fontenc}
\usepackage[utf8]{inputenc}
\usepackage{booktabs}
\usepackage{tabularx}
\usepackage{lipsum}
\usepackage{physics}
\usepackage{caption}
\usepackage{enumitem, ragged2e}
\usepackage[svgnames]{xcolor}
\usepackage{comment}
\usepackage{dsfont}
\usepackage{amsmath, amssymb, graphicx, booktabs, longtable}

\usepackage{amsthm}
\usepackage{float}
\newtheorem{thm}{Theorem}[section]
\newtheorem{lemma}[thm]{Lemma}
\newtheorem{proposition}[thm]{Proposition}
\newtheorem{definition}[thm]{Definition}

\newtheorem{remark}[thm]{Remark}

\newcommand{\R}{\mathbb{R}}

\makeatletter
\newcommand*{\compress}{\@minipagetrue}
\makeatother

\definecolor{red}{rgb}{1.0,0.0,0.0}

\definecolor{blu}{rgb}{0.0,0.0,1.0}

\definecolor{gre}{rgb}{0.03,0.50,0.03}

\definecolor{darkviolet}{rgb}{0.58, 0.0, 0.83}

\title{Optimal Policies for Environmental Assets\\ under Spatial Heterogeneity and Global Awareness}
\author{Emmanuelle Augeraud\thanks{School of Economics, University of Bordeaux, Bordeaux, France.}, Daria Ghilli\thanks{Department of Economics and Management,  University of Pavia, Pavia, italy. }, Fausto Gozzi,\thanks{Department of Economics and Finance, LUISS Guido Carli, Rome, Italy.}, Marta Leocata\thanks{Economics and Finance Department, LUISS Guido Carli, Rome, IT.}
}
\begin{document}
\maketitle

\begin{abstract}
The aim of this paper is to formulate and study a stochastic model for the management of environmental assets in a geographical context where in each place the local authorities take their policy decisions maximizing their own welfare, hence not cooperating each other.
A key feature of our model is that the welfare depends not only on the local environmental asset, but also on the global one,
making the problem much more interesting but technically much more complex to study, since strategic interaction among players arise.

We study the problem first from the $N$-players game perspective and find open and closed loop Nash equilibria in explicit form.
We also study the convergence of the $N$-players game (when $n\to +\infty$) to a suitable Mean Field Game whose unique equilibrium is exactly the limit of both the open and closed loop Nash equilibria found above, hence supporting their meaning for the game.
Then we solve explicitly the problem from the cooperative perspective of the social planner and compare its solution to the equilibria of the $N$-players game. Moreover we find the Pigouvian tax which aligns the decentralized closed loop equilibrium to the social optimum.
\end{abstract}

\vskip-0.2truecm


\vskip-1.5truecm

\tableofcontents

\vskip-1.5truecm

\section{Introduction}

Understanding the impact of human activities on the environment, and identifying effective strategies to mitigate it, represents a major challenge for contemporary scientists and policy makers. This challenge is of paramount importance at the global level, as issues such as global warming, biodiversity loss, deforestation, reduction of green areas, and pollution pose critical threats to ecosystems and human welfare.

The primary goal of the present paper is to contribute to this challenge
from the theoretical viewpoint by developing a model which incorporates the trade-off between the economic gain from exploiting the environmental assets and the benefit of preserving them and which takes into account 
new and key features with respect to the current literature.
Such new features are:
(1) the spatial heterogeneity of agents and  environment, and 
(2) the interaction between local and global environmental scales. 
Combining these two aspects requires substantial technical work but provides a richer understanding of such a complex issue. In particular it allows us to design global environmental tax policies that take into account the above two features.

Notice that (1) and (2) are deeply interconnected.
Under full information and a constant regeneration rate of environmental assets, Ayong le Kama~\cite{lekama01} and Ayong le Kama and Schubert~\cite{lekama04} showed that, in equilibrium, consumption and environmental asset grow at a constant rate, following a Balanced Growth Path (BGP). However, when heterogeneity in location and in agents’ preferences becomes relevant, agents naturally start to evaluate environmental asset on a broader, possibly global, scale. Such heterogeneity may arise from geographical or biological characteristics of different areas, or from diverse “green attitudes” among citizens and decision makers. 

\medskip

This setting leads naturally to strategic interactions among local agents, which can be formalized as an $N$-player differential game.
Now we first provide a brief account of the literature on differential games on economic-environmental topics
and then we go deeper into our model and our results on it.

\medskip

A substantial body of literature has employed differential games to address environmental decision-making problems. Among the early contributions are \cite{dockneretal, Jorgensenetal, lietal, Liuetal, long, Ploegetal}, which focus on the excessive accumulation of greenhouse gases (GHGs). In this literature, emissions are typically the control variable, and the goal is to determine the optimal emission rate that minimizes environmental damages due to GHG accumulation. Differential games have also been widely applied to decision-making problems related to watershed pollution and water taxation in agriculture \cite{biancardietal, dietal, kicsiny, misraetal, momenietal}, where water is a scarce and shared resource and the control variables are water rights or extraction levels. Another strand of literature focuses on the preservation of green areas and forests within a game-theoretic framework \cite{andres, fredj, fredj2, herran, herran2, Soest}. In these models, two main agents interact: forest owners, who exploit land for economic returns, and donors or environmentally aware players, who aim to compensate forest owners investing in the preservation of natural resources. Interestingly, in \cite{andres}, the authors enrich this framework by introducing a non-forest owner who derives utility from polluting activities.

More recently, an emerging literature has explored the application of Mean Field Games (MFGs: which are substantially limits of differential games when the number of agents goes to infinity)  to environmental economics and sustainable development. While MFGs are now well-established in economics and finance, their application to environmental problems is relatively new. For example, in \cite{kobeissietal} authors study the tragedy of the commons through an MFG approach, in \cite{lavignetankov} the decarbonization of financial markets is analyzed, and \cite{delsartoetal} pollution regulation under emission trading schemes is investigated.

\medskip

Now we go deeper into our model. It is a $N$-players differential game (together with its limit when $N \to + \infty$, i.e. its corresponding Mean Field Game), where each player can be interpreted as the planner of some local environmental assets 
(e.g. green-land/forests owners) maximizing their own welfare (i.e. a given intertemporal discounted utility function) and competing with each others.
Unlike most of the existing literature, our players exhibit a certain degree of environmental awareness: they not only care about their own consumption and local environment but also value the preservation of environmental assets across other locations. 
Hence the instantaneous utility function of the players (taken in logarithmic form) depends
on: the consumption of the player itself, the amenity value of the local environmental asset, and an index of the global environmental quality.

Concerning the state dynamics, in each location, the environmental assets regenerates at a constant rate but is negatively affected by economic activities. The value of environmental asset is also subject to both a common and an idiosyncratic source of uncertainty.


\medskip

Within this framework, we explicitly compute an open-loop Nash equilibrium (through the Pontryagin Maximum Principle) and a closed-loop Nash equilibrium (through the Dynamic Programming Approach) of the game. These equilibria may not be unique and differ by terms of order $1/N$. Although we do not establish uniqueness, we show that, as the number of players $N$ tends to infinity, both the above Nash equilibria converge to the unique equilibrium of the associated Mean Field Game. This fact confirms the relevance, for the $N$-players game, of the Nash equilibria we have computed.

Explicit computations of Nash equilibria in differential games have been done
in various contexts (see, e.g., \cite{grosset, herranrubio, lacker17, carmona2013mean}), but in our setting, as in the one of \cite{carmona2013mean}, the relationship between closed-loop and open-loop equilibria plays a crucial role. It provides strong theoretical evidence regarding equilibrium behavior in large populations and offers a confirmation of an important conjecture in the MFG literature, notably discussed by Carmona and Delarue \cite{carmona2018probabilistic1, carmona2018probabilistic2}, which states that the distinction between open-loop and closed-loop equilibria often vanishes in the limit of a large population. Our result confirms this conjecture in the context of environmental management, showing that closed-loop strategies, which capture in a sharper way the core of the strategic interactions of the finite player game, asymptotically coincide with open-loop strategies.

\medskip

Once the above Nash equilibria are computed explicitly, we analyze them 
clarifying 
under which conditions, and to what extent, the inclusion of global environmental assets in agents’ preferences affects local decision-making. 
As we expect, such inclusion improves the global environmental quality.


Then we go further since the transboundary nature of environmental problems calls for coordinated policy responses, i.e. to go beyond the decentralized decision-making across heterogeneous regions with different priorities and constraints 
and to design effective policy instruments (such as enviromental taxation) that align individual competitive behavior with collective, long-term environmental objectives.


More in detail, after the study of the decentralized (noncooperative) model, we also study the centralized (cooperative) model, i.e. the so-called "planner solution",
and compare the corresponding outcomes.
This allows to go to the main goal of this paper from the policy viewpoint: the explicit characterization of a Pigouvian tax—a corrective mechanism that internalizes the negative externalities of local consumption on the global environment—and show how it bridges the gap between the equilibrium outcomes of the strategic game and the socially optimal solution.

Concerning the literature on taxes, the theoretical foundation for corrective environmental taxation traces back to Pigou \cite{pigou} who proposed a tax to align private costs with social costs. This Pigouvian solution was  challenged by Coase in \cite{Coase}, who argued that under conditions of well-defined property rights and negligible transaction costs, bargaining between parties could achieve efficiency without central intervention. In the context of complex environmental problems like climate change, the Pigouvian approach has regained prominence. We provide selected references as examples; the relevant literature is of course much more extensive. We  refer to \cite{Wei}, which provides a fundamental rule for choosing between price (tax) and quantity instruments under uncertainty and to \cite{nordhaus}, providing a direct case for a Pigouvian carbon tax and arguing that price-type approaches such as carbon taxes have major advantages for slowing climate change. Finally we mention a cornestone in modern climate economics \cite{golosov}, deriving a precise formula for the optimal tax in a general equilibrium model including climate damages. Our analysis contributes to this tradition by explicitly characterizing a Pigouvian tax within a dynamic, multi-agent model with spatial heterogeneity.

\bigskip

\noindent Our paper is structured as follows. In Section~\ref{sec:closed_loop}, we solve the differential game using a dynamic programming approach and identify a closed-loop Nash equilibrium. Section~\ref{sec:open_loop} applies the maximum principle to characterize open-loop equilibria. Section~\ref{sec:social_planner} reformulates the problem from the perspective of a social planner. In Section~\ref{sec:tax}, we introduce a Pigouvian tax to internalize externalities and analyze how the tax rate at each location depends on the role of global environmental asset in agents’ preferences. Section~\ref{sec:uniqueness} studies the uniqueness of the equilibrium under homogeneous coefficients and the convergence of both Nash equilibria (open-loop and closed-loop) to the MFG equilibrium. Finally, Section~\ref{sec:conclusion} concludes and discusses open questions.

\section{The basic model: a game among agents indexed by locations}
 \label{sec:themodel}
We consider an economy composed of $N$ locations. In each location $i\in\left\{1,\dots,N\right\}$ the state variable at time $t$ is the local amenity value of the environmental assets  $Q_{i}(t)$, while the control variable is the percent rate $\alpha_i(t)$ at which such assets are depleted by consumption. The dynamic evolution of $Q_i(t)$ is described by the following stochastic differential equation:
\begin{equation}\label{eq:state}
\begin{cases}
dQ_{i}(t) = (\gamma_iQ_i(t)-\alpha_i(t)Q_i(t) )dt+Q_{i}(t)\mu _{i}dW^{i}\left( t\right) +Q_{i}(t)\nu _{i}dB\left(
t\right)   \\
Q_{i}\left( 0\right)=Q^i_{0}.
\end{cases}
\end{equation}
We explain now the meaning of each term.
Concerning the random terms, we assume that the dynamic of $Q_i(\cdot)$ is impacted by two exogenous phenomena: idiosyncratic shocks that are specific to each locations (typically due to geographic features) and a common noise that is due to global systemic perturbations (such as global warming).
The instantaneous volatility due to the local perturbation (denoted $\mu_{i}\geq 0)$ and the one due to global changes (denoted $\nu _{i}\geq 0)$ are specific to each location.
We call $B$ the real standard Brownian motion (SBM) associated to the common noise and $W^{i}$ the real SBM
associated to the idiosyncratic shocks in each location $i\in \left\{ 1,\dots,N\right\}$. We assume that $B,$ $W^{1},..,W^{N}$ are independent.
The parameter $\gamma _{i}>0$ is the natural regeneration rate
of the environmental assets, while $Q^i_{0}$ is the initial value of the environmental assets at location $i$.

Now we go to the constraints: as it is natural here, we assume that both the control variable and the state variable remain positive at all times. The positivity of the state variable $Q_i(\cdot)$ at all times is guaranteed by the form of the state equation, once the initial value $Q^i_{0}$ is positive. For the control variable we have to assume that $\alpha_i (t)\ge 0$ for all $t\ge 0$.

\smallskip

We are now ready to define the associated geographic game. We assume that, in each location $i=1,\dots,N$, there is a player who controls only the variable $\alpha_i$. Such $i$-th player acts as a representative agents who has to maximize her welfare, which depends on her
consumption, the level of environmental asset in her location
and an index of the global environmental asset.
We assume that this last index (denoted by $\widetilde{Q}$) is defined as the
geometric mean of the $Q_{i}$, i.e.\footnote{Observe that the term $\widetilde{Q}$  can be rewritten in term of the so-called empirical measure $\mu_N$ defined as
(here $\delta_z$ is the Dirac Delta measure centered on $z$)
$\mu_N:=\frac{1}{N}\sum_{i=1}^{N}\delta_{Q_i}$.
Hence
\[
    \widetilde{Q}=\left(\Pi_{i=1}^N Q_i\right)^{\frac{1}{N}}=e^{\frac{1}{N}\sum_{i=1}^N\ln(Q_i)
    =e^{\langle \ln(q),\mu^N\rangle}}.
\]
}
\begin{equation}
\widetilde{Q}=\left( \prod\limits_{j=1}^{N}Q_{j}\right) ^{1/N}
\label{eq:geo}
\end{equation}
More precisely, each agent aims to maximize the intertemporal utility function given by
\begin{equation}\label{eq:OCfinite}
J_i\left(\alpha_i\right) =\mathbb{E}\left[
\int_{0}^{\infty }e^{-\rho t}U_{i}\left( \alpha_{i}(t),Q_{i}(t),\widetilde{Q}(t)\right) dt
\right]
\end{equation}%
where $\rho >0$ is the discount rate and the instantaneous utility function is defined as
\[
U_{i}\left( \alpha_{i},Q_{i},\widetilde{Q}\right):=
\ln
\left[ (\alpha_{i}Q_i)(Q_i^{\theta_i})(\widetilde{Q}^{\eta_i})\right]
=
\ln \left( \alpha_{i}Q_i^{\theta_i+1}\widetilde{Q}^{\eta_i}\right)
\]
Note here that  $\theta _{i}>0$ stands for relative preference for the local environmental asset, and $\eta _{i}>0$ is the relative preference for the global environmental asset $\widetilde{Q}$. It is exactly the presence of this last term that determines the strategic interactions among the players. It is interesting to study how the presence of such term affects the decisions of players. In the case under study, where the utility has a logarithmic form, we will prove that the  effect  of such term on the growth of the local enviromental quality  is always positive.

Note also that, while the parameter $\theta_i$ is mainly related to the concrete effects of environmental asset on the welfare, the coefficient $\eta_i$ is more related to the information available to the players on the global environmental asset.\footnote{Following, e.g., Bezin \cite{bezin15} and Thaler and Sunstein \cite{thaler08} we can say that the $\eta_i$ are behavioral parameters that affect the local policies and that can be influenced e.g. increasing the media coverage about global environmental events (climate change, biodiversity losses), education.}



\subsection{Definitions of equilibria and technical lemmas}

We start with a basic definition, which is needed to introduce the notion of equilibria in our context.
\begin{definition}
Let $ \left( \mathcal{F}^i_t \right)_{t \in [0,T]}$ and $ \left( \mathcal{F}^0_t \right)_{t \in [0,T]}$ denote respectively the natural filtration generated by the Wiener processes $\mathbb{W}^i = \left( W^i_t \right)_{t \in [0,T]}$ and $\mathbb{W}^0 = \left( W^0_t \right)_{t \in [0,T]}$. Define $\left( \bar{\mathcal{F}}^i_t \right)_{t \in [0,T]}$ the complete and right-continuous augmentation of the filtration $ \left( \mathcal{F}^i_t\otimes\mathcal{F}^0_t \right)_{t \in [0,T]}$.
\noindent The set of admissible strategies for player $i$, denoted by $\mathbb{A}^i$, consists of all processes $\left( \alpha^i_t \right)_{t \in [0,T]}$ such that
\[
\int_0^T \left| \alpha^i_s \right|^2 \, ds < \infty.
\]
and $\left( \bar{\mathcal{F}}^i_t \right)_{t \in [0,T]}$- adapted.
\noindent A strategy profile for the game is a tuple $\left( \alpha^1, \dots, \alpha^N \right)$ with each $\alpha^i \in \mathbb{A}^i$. The set of admissible strategies for all players is denoted by
\[
\mathbb{A} = \mathbb{A}^1 \times \mathbb{A}^2 \times \dots \times \mathbb{A}^N.
\]
\end{definition}

We now turn to the definition of Nash equilibria, referring to Definition 5.3 and Definition 5.6 in \cite{carmona2016lectures}.

\begin{definition}\label{def:admissible_NE}
     A vector of admissible strategies $\boldsymbol{\alpha}=\left(\alpha^{ 1}, \cdots, \alpha^{ n}\right) \in \mathbb{A}$ is said to be an open-loop Nash equilibrium for the game if
$$
\forall i \in\{1, \dots, N\}, \quad \forall \alpha^i \in \mathbb{A}^i, \quad J_i\left(\boldsymbol{\alpha}\right) \geq J_i\left(\alpha^{-i}, \alpha^i\right) .
$$
We will denote such equilibrium with $\boldsymbol{\alpha}_{OL}$.
\end{definition}

\begin{definition}\label{def:CLNE}
     A vector of admissible strategies  $\boldsymbol{\alpha}=\left(\alpha^{ 1}, \cdots, \alpha^{ N}\right) \in \mathbb{A}$ is a closed-loop Nash equilibrium for the above game if, for every $i=1,\dots,N$,
     there exists a measurable function
     $\varphi^i:\R_+\times \R^N_+\to \R_+$
     such that
     \[
     \alpha^i_t=\varphi^i(t,\boldsymbol{Q}(t)), \qquad \forall t\ge 0,
     \]
     where $\boldsymbol{Q}(t)=(Q_1(t),\dots,Q_N(t))$.
     and
$$
\forall i \in\{1, \dots, N\}, \quad \forall \alpha^i \in \mathbb{A}^i, \quad J_i\left(\boldsymbol{\alpha}\right) \geq J_i\left(\alpha^{-i}, \alpha^i\right) .
$$
We will denote such equilibrium with $\boldsymbol{\alpha}_{CL}$.
\end{definition}


We now give two technical Lemmas (whose proof is in the Appendix) which allow us to rewrite the problem in a more convenient way which will be used in the next sections.
First about the rewriting of the state equation \eqref{eq:state}.
\begin{lemma}\label{lem:dynamics2new}
Equation \eqref{eq:state} can be equivalently rewritten as    \begin{eqnarray}\label{eqn:dynamicsbis}
dQ_i(t)&=&(\gamma_i Q_i(t)-\alpha_i(t)Q_i(t))dt+Q_i(t)\sqrt{\mu_i^2+\nu_i^2}dX_i(t)\\
Q_i(0)&=&Q^i_0 \quad \mbox{given}\nonumber
\end{eqnarray}
where $X^i$ is a standard Brownian motion correlated to $B$ given by
$X_i\left( t\right) =\lambda _{i}\left( \mu _{i}W^{i}\left(
t\right) +\nu _{i}B\left( t\right) \right) $ where $\lambda _{i}=\frac{1}{%
\sqrt{\mu _{i}^{2}+\nu _{i}^{2}}}.$.
\end{lemma}
\begin{proof}
    See Appendix.
\end{proof}
Second about the dynamics of the global environmental asset
in all but $i$'s location.
Let \[\widehat{Q^{-i}}=\left( \prod\limits_{j\neq i}Q_{j}\right) ^{1/N}.\]
We have
\begin{equation*}
\alpha_iQ_{i}^{\theta _{i}+1}\widetilde{Q}^{\eta _{i}}=\alpha _{i}Q_{i}^{1+\theta _{i}+\frac{\eta _{i}}{N}}\widehat{Q^{-i}}^{\eta _{i}}.
\end{equation*}
We have the following result on the dynamics of $\widehat{Q^{-i}}$.

\begin{lemma}
\label{lem:dynaQchap} The dynamics of the geometric mean of the
environmental asset in all but $i$'s location is given by
\begin{equation}\label{eqn:dynQhat}
\frac{d\widehat{Q^{-i}}}{\widehat{Q^{-i}}}=\left( \widehat{g^{-i}}+\frac{%
\widehat{\nu ^{-i}}^{2}+\widehat{\xi ^{-i}}^{2}}{2}\right) dt+\sqrt{\widehat{%
\nu ^{-i}}^{2}+\widehat{\xi ^{-i}}^{2}}d\widehat{X^{-i}(t)}
\end{equation}%
where $\sqrt{\widehat{\nu ^{-i}}^{2}+\widehat{\xi ^{-i}}^{2}}\widehat{%
X^{-i}(t)}=\widehat{\nu ^{-i}}B(t)+\frac{1}{N}\displaystyle%
\sum\limits_{j\neq i}\mu _{j}dW_{j}(t)$ is a standard Brownian motion and where%
\begin{eqnarray}
\widehat{\nu ^{-i}} &=&\frac{1}{N}\displaystyle\sum\limits_{j\neq i}\nu _{j},
\nonumber \\
\widehat{g^{-i}} &=&\frac{1}{N}\displaystyle\sum\limits_{j\neq i}\left(
\gamma _{j}-\alpha _{j}-\frac{\nu _{j}^{2}+\mu _{j}^{2}}{2}\right). \label{eqn:ghat}\\
\widehat{\xi ^{-i}} &=&\left( \displaystyle\sum\limits_{j\neq i}\frac{\mu
_{j}^{2}}{N^{2}}\right) ^{1/2}.\nonumber
\end{eqnarray}
\end{lemma}
\begin{proof}
    See Appendix.
\end{proof}
Given the above results, the problem of the agent $i$ can be reformulated in the following, technically more convenient, way: the agent chooses $\alpha_i$ that maximizes
\begin{equation*}
J_i\left( \alpha _{i},\alpha _{-i}; q,\widehat{q}\right) =\mathbb{E}\left[
\int_{0}^{\infty }e^{-\rho t}\ln\left( \alpha _{i}Q_{i}^{\left(1+\theta\right) _{i}}\widehat{Q^{-i}}%
^{\eta _{i}}\right) dt\right] ,
\end{equation*}
under the two following state equations:
\begin{eqnarray}\label{eqn:constraints}
dQ_{i}(t) &=&Q_{i}(t)\left( \left( \gamma _{i}-\alpha _{i}(t)\right) dt+\sqrt{\mu
_{i}^{2}+\nu _{i}^{2}}dX_{i}\left( t\right) \right), \quad Q_i(0)=q\\
d\widehat{Q^{-i}}(t) &=&\widehat{Q^{-i}}(t)\left( \left( \widehat{g^{-i}}+%
\frac{\widehat{\nu ^{-i}}^{2}+\widehat{\xi ^{-i}}^{2}}{2}\right) dt+\sqrt{%
\widehat{\nu ^{-i}}^{2}+\widehat{\xi ^{-i}}^{2}}d\widehat{X^{-i}(t)}\right), \quad \widehat{Q^{-i}}(0)=\widehat{q}.
\end{eqnarray}
{The value function of such problem is defined as follows:
\begin{equation}\label{eq:value_game}
V_i(q, \widehat q)=\sup_{\alpha_i}J_i\left( \alpha _{i},\alpha _{-i}; q,\widehat{q}\right)
\end{equation}
The Hamilton-Jacobi-Bellman equation associated to the above reformulation of the control problem of the player $i$ is the following:
\begin{multline}\label{eqn:HJBgame}
\rho v_i\left( q,\widehat{q}\right) -\frac{\widehat{\nu^{-i} }^{2}+\widehat{\xi^{-i}}^{2}}{2}\widehat{q}^{2}\partial_{\widehat{%
q}\widehat{q}}v_i\left( q,\widehat{q}\right) -\frac{\nu _{i}^{2}+\mu
_{i}^{2}}{2}q^{2}\partial_{qq}v_i\left( q,\widehat{q}\right)  -\widehat{q}\partial_{\widehat{q}}v_i\left( q,\widehat{q}\right)
\left( \widehat{g^{-i}}+\frac{\widehat{\nu^{-i} }^{2}+\widehat{\xi^{-i}}^{2}}{2}\right)\\-\gamma q\partial_qv_i(q,\widehat{q})-\sup_{\alpha \in \R_+}\mathcal{H}_{CV}(\alpha,\partial_qv_i,q)- \left(\nu_i \widehat{\nu^{-i}}\right) q \widehat{q} \, \partial_{q\widehat{q}}v_{i}(q, \widehat{q})=0,
\end{multline}%
where
\begin{equation}\label{eq:h_cv}
\mathcal{H}_{CV}(\alpha,\partial_qv_i,q)=
\ln(\alpha)+\left(1+\theta _{i}+\frac{\eta _{i}}{N}\right)\ln(q)+\eta_i\ln(\widehat{q})-\alpha q\partial_qv_i\left( q,\widehat{q}\right).
\end{equation}
}

\section{
A closed-loop equilibrium for the game}
\label{sec:closed_loop}

Here we provide, for every $i=1,\dots,N$, an explicit solution of \eqref{eqn:HJBgame} and exploit it to provide a closed loop equilibrium of our original game. 
The proofs of all subsequent results is given in the Appendix.
First we give a definition of solution to \eqref{eqn:HJBgame}.
\begin{definition}
    We say that a function $v_{i}$ is a classical solution of the equation \eqref{eqn:HJBgame} over $\R_{++}\times \R_{++}$ if
    $v_i\in C^2(\R_{++}\times \R_{++})$, and $v_i$ satisfies \eqref{eqn:HJBgame}.
\end{definition}
We start with following lemma. 
\begin{lemma}\label{lem:optconslingame}
Then
\begin{equation}\label{eqn:contrott2}
\alpha_i=\mbox{argmax}\left\{\ln(\alpha)-\alpha q\partial_qv_i\right\}=\left\{ \begin{aligned}
    (q\partial_qv_i)^{-1}& &  \mbox{ if } \partial_q v_i> 0\\
    \emptyset  & &\mbox{ otherwise}
\end{aligned}
\right.
\end{equation}.
\end{lemma}
\begin{proof}
    See Appendix.
\end{proof}
Now we provide an explicit solution of the HJB equation \eqref{eqn:HJBgame}.
\begin{proposition}\label{lem:vf1game}
 The function
   $$
w_i\left( q,\widehat{q}\right)  =
a\ln(q)+b\ln(\widehat q)+c,
$$
where
\begin{equation}\label{eqn:abgame}
a=\frac{\left(1+\theta _{i}+\frac{\eta _{i}}{N}\right)}{\rho}, \quad b=\frac{\eta_i}{\rho}
\end{equation}
\begin{equation}\label{eqn:cgame}
c=-\rho^{-1}\ln(\left(1+\theta _{i}+\frac{\eta _{i}}{N}\right))+\rho^{-1}\ln(\rho)-\rho^{-1}+\gamma_i \left(1+\theta _{i}+\frac{\eta _{i}}{N}\right)\rho^{-2}-\left(1+\theta _{i}+\frac{\eta _{i}}{N}\right)\rho^{-2}\frac{\nu _{i}^{2}+\mu
_{i}^{2}}{2}+\rho^{-2}\eta_i\widehat{g^{-i}}
\end{equation}
and $\widehat{g^{-i}}$ is defined in \eqref{eqn:ghat},
is a classical solution to \eqref{eqn:HJBgame}.
\end{proposition}
\begin{proof}
    See Appendix.
\end{proof}
\noindent By the proof of the previous theorem, note that we have
$$
\partial_qw_i(q, \widehat q)=\frac{a}{q}, \quad a>0.
$$
That is, the gradient of $w_i$  are positive in $(0,+\infty)$.
Therefore, due to Lemma \ref{lem:optconslingame}, Proposition \ref{lem:vf1game} we define the feedback maps as
\begin{equation}\label{eqn:Gi2}
G_i(q, \widehat q)=(q\partial_qw_{i})^{-1}=\left(\frac{\left(1+\theta _{i}+\frac{\eta _{i}}{N}\right)}{\rho}\right)^{-1}
=
\frac{\rho}{1+\theta_i+\frac{\eta_i}{N}}.
\end{equation}
We now need to know that such feedback maps are admissible, hence we have the following result.
\begin{proposition}\label{lem:feedbackgame}
For every $i=1,\dots,N$, there exists a solution to
\begin{eqnarray}\label{eqn:closedsigmagame}
dQ_i(t)&=&(\gamma_i Q_i(t)-G_i(Q_i(t), \widehat{Q^{-i}}(t)))Q_i(t)dt+Q_i(t)\mu_i dW(t)+Q_i(t)\nu_idB(t)\\
Q_i(0)&=&q^i_{0} \quad \mbox{given}
\nonumber
\end{eqnarray}
Moreover, the feedback maps $G_i$ are admissible controls.
\end{proposition}
\begin{proof}
    See Appendix.
\end{proof}

We can now pass to our main result.
\begin{thm}
\label{prop:NashN} Given $V_i$, value function defined in \eqref{eq:value_game}, we have
\begin{equation}\label{v_w}
V_i\equiv w_i, \quad i=1,\dots,N.
\end{equation}
Our game possesses a closed-loop  equilibrium given by the feedback maps
\begin{equation}\label{eqn:c3game}
G_{i}(q,\hat{q})= \rho \left(1+\theta _{i}+\frac{\eta _{i}}{N}\right)^{-1}, \qquad i=1,\dots, N
\end{equation}
Hence the resulting closed-loop equilibrium, as from Definition \ref{def:CLNE} is constant and given by
\begin{equation}\label{alpha_CL}
\alpha^i_{CL}=\rho (1+\theta_i+\eta_i/N)^{-1},
\qquad i=1,\dots, N.
\end{equation}
The associated equilibrium consumption rate at each time $t$ is given by
\begin{equation}\label{alpha_CL}
\alpha^i_{CL}Q_i(t),\qquad i=1,\dots, N.
\end{equation}
where $Q_i$ is, for $i=1,\dots,N$, the solution to
\begin{eqnarray*}
dQ_i(t)&=&(\gamma_i Q_i(t)-\rho \left(1+\theta _{i}+\frac{\eta _{i}}{N}\right)^{-1} Q_i(t)dt+Q_i(t)\mu_i dW(t)+Q_i(t)\nu_idB(t)\\
Q_i(0)&=&q^i_0
\end{eqnarray*}
i.e.
\begin{equation}\label{eqn:sol}
Q_{i}\left( t\right) =q^i_{0}e^{\left( \gamma _{i}-\rho \left(1+\theta _{i}+\frac{\eta _{i}}{N}\right)^{-1}-\frac{_{\nu
_{i}^{2}+\mu _{i}^{2}}}{2}\right) t+\nu _{i}B\left( t\right) +\mu
_{i}W^{i}\left( t\right) }.
\end{equation}
Hence, in such equilibrium, the expected growth rate in location $i$ is
$$
g_i=\gamma_i-\rho(1+\theta_i+\eta_i/N)^{-1}-\frac{\mu_i^2+\nu_i^2}{2}.
$$
and $\frac{\partial g_i}{\eta_i}>0$.
\end{thm}
\begin{proof}
    See Appendix.
\end{proof}

From the above result we deduce that, in the case under study the presence, in the utility function, of the global environment quality, has always a positive effect on the long-run dynamics of each local environmental asset.
This means that it is more convenient for the players to take into account the global environment in their objective, in order to increase the local growth rate.

Some comment on the closed loop equilibrium $\alpha^i_{CL}$. First, we observe that the optimal consumption rate decreases with respect to the local perception parameter $\theta_i$ and also with respect to the global perception parameter $\eta_i$. This means that the greater the awareness of the preservation of environmental assets, the lower the consumption rates. However, the dependence on the global variable is weak, in the sense that it becomes negligible when $N\to \infty$. This means that the model, in some sense, captures the tragedy of the commons, but the dependence on the aggregate variable remains weak. This is why we need to internalize global externalities, with a tax, see Section \ref{sec:tax}.
\section{
An open loop equilibrium for the game}
\label{sec:open_loop}
Aim of this section is to derive both open loop equilibria and then compare them with closed loop equilibria derived in the section above. The content of this section is somehow inspired by \cite{carmona2013mean}.

The $N$-player game detected by the state equations \eqref{eq:state} and by the utility function \eqref{eq:OCfinite} is here solved by looking for open-loop Nash equilibrium along Definition \ref{def:admissible_NE}.


We use, as in \cite{carmona2013mean},
the approach of Maximum Principle. Hence we start defining 
the Hamiltonians $H^i:\mathbb{R}^N\times\mathbb{R}^N\times\mathbb{R}^{N\times N}\times \mathbb{R}^N\to R$, as follows
\begin{align*}
H^i(Q,Y^i,Z_{i},\alpha)&=\langle b(Q),Y^i(t)\rangle+U_i(\alpha_i,Q_i,\tilde Q)+Trace\left(\sigma^T(Q)Z^i\right)\\
&=\sum_{k=1}^N (\gamma_kQ_k(t)-\alpha_kQ_k(t))\cdot Y^{i,k}(t)+U_i(\alpha_i,Q_i,\hat{Q})+\sum_{k=1}^N Z_{i,k,k}(t)\mu_k Q_k+\sum_{k=1}^N Z_{i,k,0}(t)\nu_kQ_k
\end{align*}
where $b,\sigma$ are respectively the drift and the diffusion of the vector $Q\in\mathbb{R}^N$ and we
consider the logarithmic utility function, therefore
\[U_i(\alpha_i,Q_i,\hat{Q})=\log(\alpha_i)+(1+\theta_i)\log(Q_i)
+\frac{\eta_i}{N}\sum_{j\neq i} \log(Q_j).\]

Similarly to what is done in \cite{carmona2016lectures}, we present the open-loop analog of the Pontryagin maximum principle.

\begin{thm}\label{thm:max_principle}
    Suppose that $\bf{\alpha}=\left(\alpha^{i}\right)_{i=1,\dots,N}$ with $\alpha^i\in \mathcal{A}_{OL}$ is an open loop Nash equilibrium, with $\left(Q^{i}\right)_{i=1,\dots,N}$ the corresponding controlled state of the system. Then, for each $i=1,\dots,N$ there exists a unique adapted solution $(Y^{i,j}_t,Z^{i,j}_t)_{t \in\mathbb{R}^+}\in \mathbb{S}^{2}\times \mathbb{H}^{2}$ of the BSDE,

    \begin{align}\label{eq:BSDE}
\begin{cases}
    dY^{i,j}(t)=\rho Y^{i,j}(t)-\left((\gamma_j-\alpha_j)Y^{i,j}(t)+\frac{(1+\theta_i)}{Q^i(t)}\delta_{i,j}+\frac{\eta_i}{N}\frac{1}{Q^j(t)}(1-\delta_{ij})+Z_{i,j,j}(t)\mu_j+Z_{i,j,0}(t)\nu_j\right)dt+\\
    \hspace{10cm}+\sum_{k=0}^N Z_{i,j,k}(t) dW^k(t)\\
        \lim_{t\to+\infty}e^{-\rho t} Q^j_t Y^{i,j}_t=0,\quad j=1,\dots,N\\
    \end{cases}
\end{align}
satisfying the coupling condition
\begin{equation}\label{eq:coupling}
-Q^i(t)Y^{i,i}(t)+\frac{1}{\alpha^i(t)}=0.
\end{equation}
Conversely, suppose that for each $i=1,\dots,N$, $\left(Q^i,\alpha^i, Y^{i,j},Z^{i,j} \right)_{j=1,\dots,N}$ is an adapted solution to the forward-backward system \eqref{eq:state}-\eqref{eq:BSDE}-\eqref{eq:coupling}. Then $\left(\alpha^i\right)_{i=1,\dots,N}$ with $\alpha^i\in \mathcal{A}_{OL}$ is an open loop Nash equilibrium.
\end{thm}
\begin{proof}
    See Appendix.
\end{proof}
    \begin{proposition}\label{prop:open_loop}
         A Nash equilibrium in open loop form is, for each $i=1,\dots,N$
        \begin{equation}\label{alpha_OL}
        \alpha^i_{OL}=\frac{\rho}{1+\theta_i}.\end{equation}
    \end{proposition}
  \begin{proof}
    See Appendix.
\end{proof}
Some comments on the open-loop equilibria $\alpha^i_{OL}$. 
First, we observe that these equilibria depend only on the local perception parameter $\theta_i$. 
As in the closed-loop case, $\alpha^i_{OL}$ decreases with $\theta_i$: the greater the awareness of local environmental preservation, the lower the consumption rate. 
However, unlike in the closed-loop setting, $\alpha^i_{OL}$ does not depend on the global perception parameter $\eta_i$. 
Indeed,
\[
\alpha^i_{CL} < \alpha^i_{OL}.
\]
This inequality can be understood by noting that closed-loop strategies reflect agents’ reactions to the evolution of the state variable $Q_i$, while open-loop equilibria do not. 
As a result, the open-loop equilibrium captures a more myopic behavior of the agents.

\section{The Social Planner solution}\label{sec:social_planner}
Objective of this section is to solve the problem presented above, from the persepective of the social planner. Therefore, we consider the following McKean Vlasov optimal control problem,
\[
J_{SP}(\boldsymbol{\alpha},Q^0)=\mathbb{E}\left[\int_0^{+\infty}e^{-\rho t}\frac1N\sum_{i=1}^N\ln(\alpha_i(t)Q_i(t)^{1+\theta_i}\widetilde{Q}^{\eta_i})\,dt\right]
\]
 $\rho>0$ is the discount rate, $Q^0\in\mathbb{R}^N$ and $Q_i$ satisfies
\begin{eqnarray}\label{eqn:dynamics}
dQ_i(t)&=&(\gamma_i Q_i(t)-\alpha_i(t)Q_i(t))dt+Q_i(t)\mu_i dW(t)+Q_i(t)\nu_idB(t)\\
Q_i(0)&=&Q_i^0.
\end{eqnarray}

and where $\boldsymbol{\alpha} \in\mathbb{A}$.
Then the value function is
\begin{equation}\label{P}
V(Q^0)=\sup_{\alpha \in \mathcal{A}(q)} J_{SP}(\alpha,Q^0).
\end{equation}

\begin{remark}
Observe that the logarithmic utility function can be rewritten as
\begin{multline*}
 \sum_{i=1}^N\ln(\alpha_i(t)Q_i(t)^{1+\theta_i}\widetilde{Q}^{\eta_i})=    \sum_{i=1}^N\left[\ln(\alpha_i(t))+(1+\theta_i)\ln(Q_i(t))+\eta_i\ln(\widetilde{Q})\right]= \\
 =\sum_{i=1}^N\ln(\alpha_i(t))+\sum_{i=1}^N(1+\theta_i)\ln(Q_i(t))+\frac{1}{N}\sum_{i=1}^N\eta_i\sum_{j=1}^N\ln(Q_j)=  \\
  =\sum_{i=1}^N\ln(\alpha_i(t))+\sum_{i=1}^N(1+\theta_i+\bar \eta)\ln(Q_i(t))
\end{multline*}
Therefore, the McKean Vlasov optimal control problem, in the logarithmic case, collapses in a standard optimal control problem.
\end{remark}

The Hamilton-Jacobi-Bellman equation associated to the McKean-Vlasov optimal control problem is
\begin{equation}\label{eqn:HJB_SP}
\rho v -\frac12\sum_{i=1}^N\left(\widehat{\nu_i }^{2}+\widehat{\xi_i }^{2}\right)q_i^{2}\partial_{q_iq_i}v-\sum_{i=1}^N\gamma_i q_i\partial_{q_i}v-\sup_{\alpha \in \R^N}\mathcal{H}_{CV}(\alpha,\partial_qv,q)=0,
\end{equation}%
where

\[
\mathcal{H}_{CV}(\alpha,\partial_qv_i,q)=
\frac1N\sum_{i=1}^N\left[\ln(\alpha_i)+(1+\theta_i+\bar{\eta})\ln(q_i)\right]-\sum_{i=1}^N\alpha_i q_i\partial_{q_i}v
\]

with $\bar{\eta}=\frac1N \sum_{i=1}^N\eta_i$.
\begin{definition}
    We say that a function $v$ is a classical solution of the \eqref{eqn:HJB_SP} over $\R^N_{++}$ if
    $V\in C^2(\R^N_{++})$, and $v$ satisfies \eqref{eqn:HJB_SP}.
\end{definition}

We will need the following lemma, which proof is omitted since it follows exactly the proof of Lemma \ref{lem:optconslingame}.
\begin{lemma}\label{lem:optconslin_SP}
\begin{equation}\label{eqn:HJB2_SP}
\alpha=\mbox{argmax}\left\{\sum_{i=1}^N\left[-\alpha_i q_i\partial_{q_i}v\right]+\frac1N\sum_{i=1}^N
\ln(\alpha_i)]\right\}
\quad\quad\text{with}\quad\alpha_i=\begin{cases}
(Nq_i\partial_{q_i}v)^{-1}\,\,\text{if}\,\,\partial_{q_i} v> 0\\
    \emptyset\hspace{1.6cm}\text{otherwise}
\end{cases}
\end{equation}
\end{lemma}
In the following theorem we find explicit solutions of the equation \eqref{eqn:HJB_SP}. Also, this proof is omitted since the proof is very similar to the proof of Proposition \ref{lem:vf1game}.
\begin{thm}\label{lem:vf_SP}
The function
    \[
    w(q)=\sum_{j=1}^N \ln{\left(a_j q_j+b_j\right)}
    \]
    with
\[
a_j=\frac{1+\theta_j+\bar{\eta}}{N\rho}, \quad
\rho b_j=-\frac12(\mu_j^2+\nu_j^2)a_j+\gamma_ja_j-1-\ln(a_j)
\]
is a classical solution to \eqref{eqn:HJB_SP}.
\end{thm}

\noindent Notice that we have
$$
\partial_{q_i}w(q)=\frac{a_i}{q_i}, \quad a_i>0.
$$
That is, the gradient of $w$  are positive in $(0,+\infty)$.
Therefore, due to Lemma \ref{lem:optconslin_SP},  and Proposition \ref{lem:vf_SP} we define the feedback maps as
\begin{equation}\label{eqn:Gi2}
G_i(q)=(\partial_{q_i}W)^{-1}, \quad i=1,\dots,N
\end{equation}

Now we prove that the functions $w$ coincide with the value functions $V$. The proof is omitted since it is the same as the proof of Theorem \ref{prop:NashN}.
\begin{thm}\label{prop:vf2game}
We have
\begin{equation}\label{eqn:uniquevsigmagame}
v=w.
\end{equation}
  Moreover, the optimal consumption is
\begin{equation}\label{eqn:c3game}
 \rho (1+\theta_i+\bar{\eta}))^{-1}Q_i(t)
\end{equation}
where $Q_i$ is the solution to
\begin{eqnarray*}
dQ_i(t)&=&(\gamma_i Q_i(t)-\rho \left(1+\theta _{i}+\frac{\eta _{i}}{N}\right)^{-1} Q_i(t)dt+Q_i(t)\mu_i dW(t)+Q_i(t)\nu_idB(t)\\
Q_i(0)&=&q_0 \quad \mbox{given}
\end{eqnarray*}
and the optimal rate of consumption is
\begin{equation}\label{alpha_SP}
\alpha^i_{SP}=\rho (1+\theta_i+\bar{\eta})^{-1}.
\end{equation}
\end{thm}
\noindent 

We observe that the social optimum consumption rate at the node $i$ depends on the local perception parameter $\theta_i$ and also on the average of the global perception parameters $\bar\eta$. As in the game theoretical framework, $\alpha^i_{SP}$ decreases with respect the parameter $\theta_i$ and $\bar\eta$. This means that the only parameter that distinguishes the optimal consumption in the node $i$ is the local perception. Indeed, the influence through the global parameter is the same in all nodes. The dependence on $\bar\eta$ captures the fact that the social planner consider the global environmental assets as a common good and therefore the global externality is fully internalized. In conclusion, we observe that
\[
\alpha^i_{OL} = \frac{\rho}{1 + \theta_i} > \alpha^i_{CL} = \frac{\rho}{1 + \theta_i + \frac{\eta_i}{N}}>\alpha^i_{SP}=\frac{\rho}{1+\theta_i+\bar{\eta}}
\]
This order between equilibria detects the progressive incorporation of environmental awareness. We start from purely local agents(open loop), then we move to adaptive but still selfish (closed loop), and we conclude with the fully cooperative case (social optimum).
\begin{proposition}\label{prop:Price_of_anarchy}
Consider the homogenous case, namely
\[\theta_i\equiv\theta,\quad \eta_i\equiv\eta,\quad \mu_i\equiv \mu,\quad \nu_i\equiv \nu, \quad i=1,\dots,N.\]
In this setting, $\alpha^i_{CL}=\alpha_{CL}=\frac{\rho}{1+\theta+\frac{\eta}{N}}$ for $i=1,\dots,N$. We denote $\boldsymbol{\alpha}_{SP}=(\alpha^i_{SP})_{i=1,\dots,N}$. Then the price of Anarchy is given by
    \[J_{SP}(\boldsymbol{\alpha}_{SP})-J^i(\alpha_{CL})=\frac{1}{\rho}\left[\frac{1+\theta+\eta}{1+\theta+\eta/N}-\log\left(\frac{1+\theta+\eta}{1+\theta+\eta/N }\right)-1\right].\]
\end{proposition}
  \begin{proof}
    See Appendix.
\end{proof}
Observe that, as $N \to \infty$, the Price of Anarchy converges to
\[
\frac{1}{\rho}\left[\frac{1+\theta+\eta}{1+\theta} - \log\!\left(\frac{1+\theta+\eta}{1+\theta}\right) - 1\right] > 0.
\]
This implies that, for a large number of players, the Price of Anarchy remains positive, meaning that the inefficiency in the competitive case persists as $N \to \infty$. 
This result is consistent with the fact that, in the closed-loop equilibria, the externality is internalized through the factor $\frac{\eta_i}{N}$. 

\noindent Moreover, we observe that the Price of Anarchy decreases with $\theta$ and increases with $\eta$. 
The decrease with respect to $\theta$ can be interpreted as follows: for fixed $\eta$, as $\theta$ increases, the optimal consumption rate (both in the game and in the social planner setting) decreases; the two rates become closer. 
Consequently, the social planner has less room for improvement, and the Price of Anarchy decreases with $\theta$. 

\noindent On the other hand, the increase with respect to $\eta$ can be explained by noting that a higher value of $\eta$ lowers the social optimum, leading the two optimal consumption rates to diverge. Therefore the Price of Anarchy increases with $\eta$.

\section{Introducing a Tax on Consumption}
\label{sec:tax}
In this section, we extend the decentralized framework by incorporating a tax on consumption. This instrument is intended to correct the misalignment between private incentives and the social cost of environmental degradation. Formally, we assume that each local authority sets a constant consumption tax rate $\tau_i \in [0,1)$ on the intensity of polluting activities. The tax is levied on the consumption-intensive component of consumption, namely $\alpha_i Q_i$, and enters the utility function as a multiplicative distortion.

The agent’s objective thus becomes:
\begin{equation}
J_i (\alpha_i) = \mathbb{E} \int_0^\infty e^{-\rho t} \left[\ln \left( (\alpha_i(t) Q_i(t)) Q_i(t)^{\theta_i} \widetilde{Q}(t)^{\eta_i}\right) -\tau_i\alpha_i\right]dt,
\end{equation}
where $\widetilde{Q}(t) = \left( \prod_{j=1}^n Q_j(t) \right)^{1/n}$ denotes the global environmental asset perceived as a geometric mean of all local qualities.

The introduction of the tax reduces the marginal utility of consumption and thus induces a downward adjustment in the optimal consumption intensity. Importantly, the tax does not directly alter the dynamics of environmental asset, which remain governed by Equation \eqref{eq:state}.

This setting allows us to explore how fiscal instruments can influence local behavior in a stochastic environment characterized by both idiosyncratic and systemic environmental shocks. It also enables a comparison of decentralized equilibria with and without taxation, and an assessment of the extent to which taxation improves individual and collective environmental outcomes.
Let us rewrite the HJB associated to the $N$-players game:
\begin{multline}
\rho v_i\left( q,\widehat{q}\right) -\frac{\widehat{\nu^{-i} }^{2}+\widehat{\xi^{-i}}^{2}}{2}\widehat{q}^{2}\partial_{\widehat{%
q}\widehat{q}}v_i\left( q,\widehat{q}\right) -\frac{\nu _{i}^{2}+\mu
_{i}^{2}}{2}q^{2}\partial_{qq}v_i\left( q,\widehat{q}\right)  -\widehat{q}\partial_{\widehat{q}}v_i\left( q,\widehat{q}\right)
\left( \widehat{g^{-i}}+\frac{\widehat{\nu^{-i} }^{2}+\widehat{\xi^{-i}}^{2}}{2}\right)\\-\gamma q\partial_qv_i(q,\widehat{q})-\sup_{\alpha \in \R_+}\mathcal{H}_{CV}(\alpha,\partial_qv_i,q)- \left(\nu_i \widehat{\nu^{-i}}\right) q \widehat{q} \, \partial_{q\widehat{q}}v_{i}(q, \widehat{q})=0,
\end{multline}%
where
$$
\mathcal{H}_{CV}(\alpha,\partial_qv_i,q)=
\ln(\alpha)+\left(1+\theta _{i}+\frac{\eta _{i}}{N}\right)\ln(q)+\eta_i\ln(\widehat{q})-\alpha q\partial_qv_i\left( q,\widehat{q}\right)
$$
By adapting results presented in Section \ref{sec:closed_loop}, we get the following result.
\begin{lemma}[Decentralized Optimum]
The decentralized closed loop equilibrium is given by:
\begin{align}
\alpha_i^{CL}(\tau_i) &= \frac{1}{\tau_i+\frac{(1+\theta_i+\frac{\eta_i}{N})}{\rho}}\label{eq:cl}
\end{align}
\end{lemma}

By imposing the equality between consumption rates $\alpha_i^{CL} (\tau^i)= \alpha_i^{SP}$, we have \eqref{eq:tau}, we derive the optimal tax rule.
\begin{proposition}[Pigouvian Tax]
The optimal tax aligning decentralized and social optimum is:
\begin{align}\label{eq:tau}
\tau_i &= \frac{\bar{\eta}-\frac{\eta_i}{N}}{\rho}.
\end{align}
\end{proposition}
The Pigouvian tax of player $i$ internalizes the social value of the global environmental asset by incorporating the average awareness of all other agents on global environmental assets. 
Note that $\tau_i$ does not depend directly on the local parameter $\theta_i$, but since $\frac{\partial \tau_i}{\partial \eta_j} > 0$, an increase in the global environmental awareness of any other player $j$ raises the tax for agent $i$. 
With respect to the discount rate, the tax is decreasing in $\rho$; thus, the more forward-looking (i.e., the lower the discount rate) the planner is, the higher the optimal tax. 
This reflects intergenerational fairness: societies that assign greater value to future welfare tend to impose stricter environmental taxation in the present.


\section{On uniqueness of equilibria for large $N$: the help of Mean Field Games  
}\label{sec:uniqueness}
In the following section we investigate a way of characterizing the closed loop equilibria that we studied in the previous sections. We underline that we have no uniqueness result for such equilibria.
Up to our knowledge, uniqueness of equilibria are typically proven by studying uniqueness of the HJB associated. This result is not trivial in our case due to the degeneracy on the boundary of the second order term, the constraint and the non Lipschitz coefficients.

 In the present section we prove uniqueness of the Mean Field Game associated and a convergence result of the equilibrium in the N-players game to the solution of the Mean Field Game. Then, since the solution to the Mean Field Game is unique, our equilibrium will exacly be the one which is selected by the Mean Field Game in the limit as the number of the agents grows. This result does not give uniqueness of the equilibria in the N-players game. However it is a characterization of the closed loop equilibria which were studied in the previous sections.
\subsection{The Mean Field Game associated to our game}

We focus on the case in which the agents are all identical. Therefore we are in the homogenous setting, with 
\[\theta_i\equiv\theta,\quad \eta_i\equiv\eta,\quad \mu_i\equiv \mu,\quad \nu_i\equiv \nu, \quad i=1,\dots,N.\]
For simplicity of exposition we focus on the case in which there is no common noise. Then the   dynamics of a reprensentative player is
\begin{equation}
\begin{cases}
dQ(t) = (\gamma Q(t)-\alpha(t)Q(t) )dt+\mu Q(t) dW\left( t\right)  \\
Q\left( 0\right)=Q_{0}.
\end{cases}
\end{equation}
We define the space
$$
 \mathcal{P}_l(\R_+)=\left\{m : m \mbox{ is a probability measure on } \R_+ \mbox{ and } \int_0^{+\infty} \ln(Q)m(t,dQ)<+\infty\right\}.
$$
where $\R_+=(0,+\infty)$.
Let  $m \in \mathcal{P}_l(\R_+)$. Then each agent aims to maximize the utility function given by
\begin{equation}\label{eq:OCfinite}
J\left(\alpha\right) =\mathbb{E}\left[
\int_{0}^{\infty }e^{-\rho t}\left(\ln ( \alpha )+(\theta+q)\ln(Q)+ \eta \langle \ln(Q), m(t)\rangle\right) dt
\right]
\end{equation}%
where $\rho >0$ is the discount rate and
$$
\langle \ln(Q), m(t)\rangle=\int_0^{+\infty} \ln(Q)m(t,dQ).$$

\begin{remark}
    Note that as proved by \eqref{eqn:sol} it is straightforward to see that the process $Q(\cdot)$ always remains positive. This is consistent with our expectation, since we expect enviromental quality to be positive. Then the state space of the probability measure $m$ is $(0,+\infty)$.
\end{remark}
\begin{remark}
    Note that by the law of large numbers we have that
    $$
    m^N_t \to m(t) \mbox{ as } N \to + \infty \mbox{ in } \mathcal{P}_l(\R_+)
    $$
    where $m(t)=\mathcal{L}(Q_t)$ and  $\mathbb{E}[Q_t]=\int_0^{+\infty}qdm(q)$ and $\mathcal{P}_l(\R_+)$ is equipped with the weak topology. Then one can expect that
    \begin{equation}\label{eq:aggregation_term_rewriting}
    \ln(\tilde Q)=\frac1n \sum_{i=1}^N\ln(Q_i)=\langle \ln(q),m^N_t\rangle \to \langle \ln(Q), m(t)\rangle
\end{equation}
\end{remark}

We write the associated Mean Field Game system:
\begin{equation}\label{MFGeq}
\begin{cases}
\rho u(q)-\frac{\mu^2}{2}q^2\partial_{qq}^2 u(q)=H(q,\alpha, \partial_q u)+\gamma q \partial_q u(q)+(1+\theta)\ln(q)+\eta \langle \ln(q),m\rangle\\
 m_t(t,q)-\frac{\mu^2}{2}q^2\partial_{qq}^2m(q)+\mbox{ div }(H_p(\alpha,q, \partial_q u)m(t,q))=0
 \end{cases}
\end{equation}
with $m(0,q)=m_0$ and
where
\begin{equation}\label{eqn:H}
H(q,\alpha,p)=\sup_{\alpha}\left\{ \ln(\alpha)-\alpha q p\right\}.
\end{equation}

\subsection{The convergence of the N-players game to the  Mean Field Game}
In the present paragraph we prove the convergence of Nash equilibria for the $N$-players game to the solution of the associated Mean Field Game \eqref{MFGeq}.

We recall that under some assumptions and in some specific cases, it has been proved that any solution to a Mean Field Game corresponds to an \(\varepsilon\)-Nash equilibrium of the associated \(N\)-player game (see \cite{carmona2018probabilistic1}, Part II, Chapter 6, Section 6.1). Moreover, in specific cases, it is possible to prove the convergence of Nash equilibria from the finite \(N\)-player game to the solution of the corresponding Mean Field Game (see \cite{carmona2018probabilistic2}, Part II, Chapter 6, Sections 6.2 and 6.3, and \cite{carmona2018probabilistic2}, Chapter 8, Section 8.2).

Later on,  Lions proved  that the solutions of the Mean Field Game are just the trajectories of a new infinite-dimensional (even if the state space has a finite dimension) PDE in the space of measures, which is called \emph{Master Equation} (ME hereafter). In other words, the ME is a partial differential equation that governs the evolution of the value function and the distribution of agents in a large population of interacting decision-makers. It encapsulates the equilibrium dynamics by linking the optimal control of individual agents with the overall distribution of the population.
Thus, the ME is a fundamental object to study to understand the properties of the discrete model's convergence to the continuous macroscopic Mean Field Game.

We will exploit the ME in order to prove the convergence result, namely Theorem \ref{thm:convme}.
Morover we will prove a uniqueness result for the solution to the ME and from this result deduce the uniqueness of the associated Mean Field Game.

\subsubsection{Some basic notion}
 To derive the ME, we fix an initial condition $(t_0,m_0)\in [0,T)\times\mathcal P_l(\R_+)$, we consider the solution $(u,m)$ of \eqref{MFGeq} and we define a function
$
U:[0,T]\times \R_+\times\mathcal P_l(\R_+)\to\R
$
as
\[
U(t,x,m_0) := u(t,x)
\]
where $u$ is the solution to \eqref{MFGeq} with initial condition $m_0$.

To compute, at least formally, the equation satisfied by $U$ (hence, the ME), we need to give a suitable definition of the derivative of $U$ with respect to the measure variable.
\begin{definition}
Let $U:\mathcal P_l(\R_+)\to\R$. We say that $U$ is $\mathcal C^1$ in the measure variable if there exists a map $K:\mathcal P_l(\R_+)\times \R_+\to\R$ such that, for all $m_1$, $m_2\in\mathcal P_l(\R_+)$, it holds
\begin{equation}\label{def:dmV}
\lim\limits_{s\to0^+}\frac{U(m_1+s(m_2-m_1))-U(m_1)}s=\int_0^{+\infty} K(m_1,\xi)\,d(m_2-m_1)(\xi)\,.
\end{equation}
We call $\frac{\partial U}{\partial m}(m,\xi)$ the unique $K$ satisfying \eqref{def:dmV} and
\begin{equation}\label{eq:normalizing}
\int_0^{+\infty} K(m,\xi)\,dm(\xi)=0\qquad\forall m\in\mathcal P_l(\R_+)\,.
\end{equation}
Moreover, if $\frac{\partial U}{\partial m}(m,\cdot)$ is $C^1$ in the space variable, we define $D_mU:\mathcal P_l(\R_+)\times \R_+\to \R$ the intrinsic derivative of $U$ as
$$
D_mU(m,\xi)=D_\xi\frac{\partial U}{\partial m}(m,\xi)\,.
$$
\end{definition}
\subsubsection{Deriving the Master Equation}
Now let us write the Nash system for the $N$-players game and see how we can derive the ME from such system. The Nash system associated to the $N$-players game in the case of no presence of common noise is the following:
\begin{equation}\label{eqn:nash}
\rho v_i-\frac{\mu^2}{2}\sum\limits_{i=1}^Nq_i^2\partial^2_{q_i q_i}v_i-(\gamma q_i-\alpha^*q_i) \partial_{q_i} v_i
-\sum\limits_{j \neq i}( \gamma q_j-\alpha^* q_j) \partial_{q_j} v_i= \ln(\alpha^*)+\left(1+\theta\right) \ln(q_i)+\eta \ln(\widehat q)
\end{equation}
We can derive, informally, the ME from the previous equation. The intuition is that we can make the following replacements
\begin{equation}\label{eqn:firstrepla}
\partial_{q_i}v_i\rightarrow \partial_q U, \quad \partial_{q_iq_i}v^i\rightarrow \partial_{qq}U
\end{equation}
and for $j \neq i$
\begin{equation}\label{eqn:secondrepla}
\partial_{q_j q_j}v_i\rightarrow\partial_vD_\mu U, \quad \partial_{q_j}v_i\rightarrow D_\mu U
\end{equation}
By \eqref{eqn:secondrepla} this way we can do the following replacement
$$
\sum_{j \neq i}q_j^2\partial_{q_jq_j}^2v_i\rightarrow  \int_0^{+\infty}(q^2 \partial_v D_\mu U) \, d\mu(v)
$$
$$
 \sum_{j \neq i}(\gamma q_j-\alpha^*qj)\partial_{q_j}v_i\rightarrow \int_0^{+\infty} (\gamma q-\alpha^*q)  D_\mu U \, d\mu(q)
$$

Note that the other terms in \eqref{eqn:nash} can be found very easily through \eqref{eqn:firstrepla}.

We remark that the above argument is not rigourous. We are just interested in giving  the idea of the procedure usually done to find, informally, the ME and the flavour on how the nonlocal terms appear in the ME.

Then, the ME is
\begin{multline}\label{eqn:me}
\rho U- \frac{\mu^2}{2} q^2\partial_{qq}^2 U- \frac{\mu^2}{2} \int_0^{+\infty}(q^2 \partial_v D_\mu U) \, d\mu(v)-(\gamma q-\alpha^*q)  \partial_q U - \int_0^{+\infty} (\gamma q-\alpha^*q)  D_\mu U \, d\mu(q)
= \\=\ln(\alpha^*)+\left(1+\theta\right) \ln(q)+\eta \langle \ln(q),\mu\rangle
\end{multline}

\subsubsection{Explicit solution of the Master Equation}

We solve the ME by a guess on the solution. Let us make the following ansatz
\begin{equation}\label{eqn:Ulog}
U(q,m) = a \ln q + b \langle \log(q), m \rangle +c
\end{equation}

The spatial derivatives are
\[
\partial_q U = \frac{a}{q}, \quad \partial_{qq}^2 U = -\frac{a}{q^2},
\]

the Lions derivative is
\[
D_\mu U = \partial_q (\partial_\mu U) = \partial_q (b\log q) = \frac{b}{q}
\]

and the velocity derivative is
\[
\partial_q D_\mu U = -\frac{b}{q^2}.
\]

Substituting our derivatives into the ME \eqref{eqn:me} we get
\[
-\rho a \ln(q)-\rho b \langle \log(q),m\rangle-\rho c+a\gamma -a\alpha^* + \gamma b -\alpha^*b -\frac{b\mu^2}{2}- \frac{a \mu^2}{2 } + \ln(\alpha^*)+\left(1+\theta\right) \ln(q)+\eta  \langle \log(q),m\rangle = 0
\]
Since \begin{equation}\label{eq:mfg_control}
\alpha^*=\frac{1}{a},
\end{equation}
\[
-\rho a \ln(q)-\rho b \langle \log(q),m\rangle-\rho c+a\gamma -1 + \gamma b -\frac{1}{a}b -\frac{b\mu^2}{2}- \frac{a \mu^2}{2 } + \ln(\frac{1}{a})+\left(1+\theta\right) \ln(q)+\eta  \langle \log(q),m\rangle = 0
\]

Then
$$
\rho a=\left(1+\theta\right), \quad \rho b=\eta
$$
$$
-\rho c+\frac{\left(1+\theta\right)}{\rho}\gamma-1+\gamma \frac{\eta}{\rho}
-\left(1+\theta\right)^{-1}\eta-\frac{\left(1+\theta\right)\mu^2}{2\rho}-\frac{\eta\mu^2}{2\rho}+\ln(\rho)-\ln( \left(1+\theta\right))=0
$$
that is
\[
 c=-\rho^{-1}\ln( \left(1+\theta\right))+\rho^{-1}\ln(\rho)-\rho^{-1}+\gamma\left(1+\theta\right)\rho^{-2}-\rho^{-2}\frac{\left(1+\theta\right)\mu^2}{2}+\gamma \eta \rho^{-2}
-\left(1+\theta\right)^{-1}\rho^{-1}\eta-\rho^{-2}\frac{\eta\mu^2}{2}
\]

Since the coefficients $a,b,c$ are uniquely determined, we  just proved that there exists a unique solution to the ME within the class of function of the type \eqref{eqn:Ulog}.

\subsubsection{The convergence  and the uniqueness result}
 Take $w$ solution of the nash system
 \begin{equation}\label{vni_explicit}
w_i\left( q,\widehat{q}\right)  =
a\ln(q)+b\ln(\widehat q)+\tilde c,
\end{equation}
where
\begin{equation}\label{eqn:abgame}
a=\frac{\left(1+\theta\right)}{\rho}, \quad b=\frac{\eta}{\rho}
\end{equation}
\begin{equation}\label{eqn:cgame}
\tilde c=-\rho^{-1}\ln(\left(1+\theta\right))+\rho^{-1}\ln(\rho)-\rho^{-1}+\gamma \left(1+\theta\right) \rho^{-2}-\rho^{-2}\frac{\left(1+\theta\right)\mu^{2}}{2}+\rho^{-2}\eta\widehat{g^{-i}}
\end{equation}

The convergence result is proven in the following theorem.
\begin{thm}\label{thm:convme}
 If for all $i=1,\dots,N$ we define the function $w^{N,i}:R\times\mathcal P_l(\R_+)\to\R$ in the following way:
$$
w^{N,i}(q_i,m):=\int\limits_{\R^{N-1}}w_i(q,\hat q) \prod\limits_{j\neq i}m(dq_j)\,,
$$
then we have, for $N\gg1$ and for a certain constant $C>0$,
\begin{equation}\label{convergenza2}
\sup\limits_{\substack{1\le i\le N}}\int_{\R_+} \big|w^{N,i}(q,m)-U(q,m)\big|\,m(dq)\le \frac{C}{N}\left(1+ \left| \int_{\R_+}\ln(q) m(dq)\right|\right) \qquad\forall\,m\in\mathcal P_l(\R_+)\,,
\end{equation}
where $C$ is a constant depending on the data.
\end{thm}
\begin{proof}
    See Appendix.
\end{proof}
\begin{remark}
    The convergence is local in $m$, due to the term $\int_0^{+\infty} \ln(q)m(dq)$. Note that this term is finite by assumption  since $m \in \mathcal{P}_l(\R_+)$.
\end{remark}

As stated in the beginning of the section and as prove in the following Theorem \ref{thm:unique} we have uniqueness of the solution to the Mean Field Game \eqref{MFGeq}. The result stems from the Lasry-Lions monotonicity condition, typically ensuring uniqueness in Mean Field Game systems and which is satisfied in our case.
Recall that the interaction is through the mean of log environmental asset:

\[
F(m) = \int_0^{+\infty} \ln(q) \, m(t, dq),
\]

In the Lasry-Lions sense, an interaction term $F(m)$ is monotone if for any two probability measures $m_1, m_2$:\[
\int \big( F(m_1) - F(m_2) \big) \, d(m_1 - m_2) \geq 0.
\]

In our case $F(m) = \int_0^{+\infty} \ln(q) \, m(dq)$ is a scalar, so we have

\[
\int_0^{+\infty}\big( F(m_1) - F(m_2) \big) \, d(m_1 - m_2)
= \big( F(m_1) - F(m_2) \big) \cdot \int_0^{+\infty} d(m_1 - m_2).
\]

But $\int_0^{+\infty} d(m_1 - m_2) = 1 - 1 = 0$.
So:

\begin{equation}\label{eqn:llmon}
\int \big( F(m_1) - F(m_2) \big) \, d(m_1 - m_2) = 0 \quad \text{for all } m_1, m_2.
\end{equation}
Once we have \eqref{eqn:llmon}, we can apply in Theorem 3.29 of \cite{carmona2013mean} to deduce uniqueness. Note that the assumptions of the above mentioned theorem are satisfied in our case by \eqref{eqn:llmon} and since the Hamiltonian in \eqref{eqn:H} is stricly concave in $\alpha$. This result, together with the convergence proved in the previous theorem, characterizes the closed loop equilibrium for the $N$-players game. In other words, the closed loop equilibrium that we found is exactly the one which converges to the solution to the Mean Field Game associated when the number of players tends to infinity.
\begin{thm}\label{thm:unique}
    Under our standing assumptions and under \eqref{eqn:llmon}, there is a unique solution to the Mean Field Game system \eqref{MFGeq}.
\end{thm}

\noindent From Theorems \ref{thm:convme} and \ref{thm:unique}, we deduce that the value functions of the $N$-player differential game, $w^{N,i}$, converge to $U(q,m)$, the unique solution of the Master equation~\eqref{eqn:me}. 
By combining~\eqref{eq:mfg_control} and~\eqref{eqn:abgame}, it follows that the unique Mean Field Nash equilibrium is given by
\[
\alpha^* = \frac{\rho}{1+\theta}.
\]
This result enables us to identify equilibria in the $N$-player game. 
Indeed, at the finite-player level, the uniqueness of equilibria is not guaranteed. 
However, in this homogeneous setting
we have
\[
\alpha^i_{CL} \to \alpha^*, \quad N \to \infty,
\]
and
\[
\alpha^i_{OL} \equiv \alpha^*, \quad i = 1,\dots,N.
\]
Hence, among the possibly multiple equilibria of the finite-player game, $\alpha^i_{CL}$ and $\alpha^i_{OL}$ can be regarded as meaningful equilibria, since they converge in the limit to the unique Mean Field Nash equilibrium.

\section{Conclusion}\label{sec:conclusion}
In this paper
we first studied an $N$-players game where agents are spatially heterogenous  but take into account the global environmental quality in their utility.
Using also the Mean Field Game theory, we investigate the Nash equilibria (open-loop and closed-loop) of such games and, in particular, the positive effect of the global awareness of the agents on the environmental quality in the equilibria.
Then we solve the problem on the perspective of the social planner and find (as expected) that the optimum of the social planner is always smaller than both the closed loop and open loop $N$-players optima. We then extend the decentralized framework by introducing a  tax on consumption and find the Pigouvian tax, aligning decentralized and social optima. 

\noindent 
For future developments of this line of research, we aim to study more general models of environmental assets by addressing two main aspects. 
First, we plan to incorporate stronger spatial interactions among agents. 
Second, we intend to explore nonlinear growth dynamics, such as logistic or Gompertz dynamics, which are particularly suitable for modeling other environmental quantities, such as biodiversity.

\newpage

\section*{Appendix}

\subsection*{Proof of Lemma \ref{lem:dynamics2new}}

Let us consider the biodiversity dynamics, given by Equation \ref%
{eqn:dynamics}:
\begin{equation*}
dQ_{i}(t)=\left( \gamma _{i}Q_{i}(t)-\alpha_{i}Q_i(t)\right) dt+Q_{i}(t)\mu _{i}dW^{i}\left(
t\right) +Q_{i}(t)\nu _{i}dB\left( t\right)
\end{equation*}

Let us define $X_i\left( t\right) =\lambda _{i}\left( \mu _{i}W^{i}\left(
t\right) +\nu _{i}B\left( t\right) \right) $ where $\lambda _{i}=\frac{1}{%
\sqrt{\mu _{i}^{2}+\nu _{i}^{2}}}.$ We are now going to prove that $X_i\left(
t\right) $ is a Brownian motion. First, as $W_{i}$ and $B$ are continuous, $%
X_i $ is continuous and moreover $X_i(0)=0.$ As $W_{i}$ and $B$ are $\mathcal{N}%
\left( 0,t\right) ,$ then $\mathbb{E}\left[ X_i\left( t\right) \right]
=\lambda _{i}\left( \mu _{i}\mathbb{E}\left[ W_{i}\left( t\right) \right]
+\nu _{i}\mathbb{E}\left[ B\left( t\right) \right] \right) =0$ and%
\begin{equation*}
\mathbb{E}\left[ X_i\left( t\right) ^{2}\right] =\lambda _{i}^{2}\left( \mu
_{i}^{2}\mathbb{E}\left[ W_i{2}\left( t\right) \right] +\nu _{i}^{2}\mathbb{E%
}\left[ B^{2}\left( t\right) \right] \right) +\lambda _{i}\mu _{i}\nu _{i}%
\mathbb{E}\left[ B\left( t\right) W_{i}\left( t\right) \right]
\end{equation*}

As $W_{i}$ and $B$ are independent, $\mathbb{E}\left[ B\left( t\right)
W_{i}\left( t\right) \right] =0$, thus $\mathbb{E}\left[ X_i\left( t\right)
^{2}\right] =\lambda _{i}^{2}\left( \mu _{i}^{2}+\nu _{i}^{2}\right) t.$ As
we have set that $\lambda _{i}=\frac{1}{\sqrt{\mu _{i}^{2}+\nu _{i}^{2}}},$ $%
\mathbb{E}\left[ X_i\left( t\right) ^{2}\right] =t.$ It can be easily checked
that $X_i(t)-X_i(s)\sim \mathcal{N}\left( 0,t-s\right) ,$ $0\leq s\leq t.$

\subsection*{Proof of Lemma \ref{lem:dynaQchap}}

\begin{proof}
As there is no risk of mistaking, to ease the reading of the proof we will
write $\widehat{x}$ instead of $\widehat{x^{-i}}$ when variable or parameter
$x$ is concerned.

$\widehat{Q}$ is a state variable defined as

\begin{equation*}
\widehat{Q}=\left( \prod\limits_{j\neq i}Q_{j}\right) ^{1/n}
\end{equation*}%
Thus $\ln \left( \widehat{Q}\right) =\frac{1}{n}\displaystyle%
\sum\limits_{j\neq i}\ln Q_{j}$. We let $Y_{i}=\ln \left( Q_{i}\right) .$%
Using Ito formula yields
\begin{equation*}
dY_{j}(t)=\frac{dQ_{j}(t)}{Q_{j}(t)}-\frac{\nu _{j}^{2}+\mu _{j}^{2}}{2}dt.
\end{equation*}%
As
\begin{equation*}
\frac{dQ_{i}(t)}{Q_{i}(t)}=\left( \gamma _{i}-\alpha _{i}\right) dt+\mu
_{i}dW_{i}\left( t\right) +\nu _{i}dB\left( t\right) ,
\end{equation*}%
then
\begin{equation*}
dY_{j}(t)=\left( \gamma _{j}-\alpha _{j}-\frac{\nu _{j}^{2}+\mu _{j}^{2}}{2}%
\right) dt+\nu _{j}dB(t)+\mu _{j}dW_{j}(t).
\end{equation*}

Let us now define $\widehat{Y}=\frac{1}{n}\sum\limits_{j\neq i}Y_{j}$ and
let
\begin{eqnarray*}
\widehat{g} &=&\frac{1}{n}\displaystyle\sum\limits_{j\neq i}\left(
\gamma _{j}-\alpha _{j}-\frac{\nu _{j}^{2}+\mu _{j}^{2}}{2}\right) , \\
\widehat{\nu ^{-i}} &=&\frac{1}{n}\displaystyle\sum\limits_{j\neq i}\nu _{j},
\end{eqnarray*}%
then the dynamics of $\widehat{Y}$ is then given by
\begin{equation*}
d\widehat{Y}(t)=\widehat{g}dt+\widehat{\nu }dB(t)+\frac{1}{n}\displaystyle%
\sum\limits_{j\neq i}\mu _{j}dW_{j}(t).
\end{equation*}

We are then going to prove that there exists a Brownian motion, thereafter
denoted $\widehat{Z},$ such that
\begin{equation*}
\left( \displaystyle\sum\limits_{j\neq i}\frac{\mu _{j}^{2}}{n^{2}}\right)
^{1/2}d\widehat{Z}_{t}=\frac{1}{n}\displaystyle\sum\limits_{j\neq i}\mu
_{j}dW_{j}(t).
\end{equation*}

Indeed, letting $\left( \displaystyle\sum\limits_{j\neq i}\frac{\mu _{j}^{2}%
}{n^{2}}\right) ^{-1/2}=\chi ,$ as $\left( W_{j}(t)\right) _{j\neq i}$ are
continuous, $\widehat{Z}_{t}$ is continuous, and moreover $\widehat{Z}%
(0)=0. $

As $\left( W_{j}(t)\right) _{j\neq i}$ are independant Brownian, $\mathbb{E}%
\left[ \widehat{Z}(t)\right] =\chi \frac{1}{n}\displaystyle%
\sum\limits_{j\neq i}\mu _{j}\mathbb{E}\left[ W_{j}(t)\right] =0.$
Moreover, variance is given as follows%
\begin{eqnarray*}
Var\left( \widehat{Z}(t)\right) &=&\mathbb{E}\left[ \widehat{Z}(t)^{2}%
\right] \\
&=&\frac{\chi ^{2}}{n^{2}}\left( \displaystyle\sum\limits_{j\neq i}\mu
_{j}^{2}\mathbb{E}\left[ \left( W_{t}^{j}\right) ^{2}\right] +\sum\limits
_{\substack{ j\neq k  \\ j\neq i,k\neq i}}\mu _{j}\mu _{k}\mathbb{E}\left[
W{j}(t)W_{k}(t)\right] \right) .
\end{eqnarray*}%
Due to independance of $\left( W_{t}^{j}\right) _{j\neq i}$, $\mathbb{E}%
\left[ W_{j}(t)W_{k}(t)\right] =0$ and thus $Var\left( \widehat{Z}%
_{t}\right) =t.$

\bigskip As a consequence,
\begin{equation*}
d\widehat{Y}(t)=\widehat{g}dt+\widehat{\nu }dB(t)+\widehat{\xi }d\widehat{Z%
}(t),
\end{equation*}%
where $\widehat{\xi }=\left( \displaystyle\sum\limits_{j\neq i}\frac{\mu
_{j}^{2}}{n^{2}}\right) ^{1/2}.$ As $B\left( t\right) $ and $Z\left(
t\right) $ are independent, similarly as previously, it can be obtained that
\begin{equation*}
d\widehat{Y}(t)=\widehat{g}dt+\sqrt{\widehat{\nu }^{2}+\widehat{\xi }^{2}}d%
\widehat{X}(t)
\end{equation*}%
where $\sqrt{\widehat{\nu }^{2}+\widehat{\xi }^{2}}$ $\widehat{X}=\widehat{%
\nu }B(t)+\widehat{\xi }\widehat{Z}(t)$ is a standard Brownian motion.

As $\widehat{Y}=\ln \left( \widehat{Q}\right) ,$ it follows that $\frac{d%
\widehat{Q}(t)}{\widehat{Q}(t)}=d\widehat{Y}(t)+\frac{\widehat{\nu }^{2}+\widehat{\xi
}^{2}}{2}dt,$ thus:

\begin{equation*}
\frac{d\widehat{Q}(t)}{\widehat{Q}(t)}=\left( \widehat{g}+\frac{\widehat{\nu }%
^{2}+\widehat{\xi }^{2}}{2}\right) dt+\sqrt{\widehat{\nu }^{2}+\widehat{\xi }%
^{2}}d\widehat{X}(t).
\end{equation*}
\end{proof}


\subsection*{Proof of Lemma \ref{lem:optconslingame}}

\begin{proof}
We derive the current value hamiltonian $\mathcal{H}_{CV}$ with respect to $\alpha$ and we get
$$
-q\partial_qv_{i}+\frac{1}{\alpha}
$$
Again, if $\partial_q v_{i}<0$ then the above expression is always positive and the expression in brackets increasing in $\alpha$. Otherwise we impose
$$
-q\partial_qv_{i}+\frac{1}{\alpha}=0
$$
and we find
\begin{equation}\label{eqn:c11}
\alpha=(q\partial_qv_{i})^{-1}.
\end{equation}
Computing the second derivative w.r.t. $\alpha$ of the $\mathcal{H}_{CV}$ and computing it in \eqref{eqn:c11} shows that \eqref{eqn:c11} is indeed a maximum.
    \end{proof}

\subsection*{Proof of Proposition \ref{lem:vf1game}}
\begin{proof}
First we  compute the partial derivatives of the candidate function. Our candidate value function is:
\[
w_i(q, \hat{q}) = a \ln(q) + b \ln(\hat{q}) + c
\]
The first derivative with respect to $q$ is:
\[
\partial_q w_i = \frac{a}{q}, \quad q \partial_q w_i = a
\]
The second derivative with respect to $q$ is:
\[
\partial_{qq} w_i = -\frac{a}{q^2}, \quad q^2 \partial_{qq} w_i = -a
\]
The first derivative with respect to $\hat{q}$ is:
\[
\partial_{\hat{q}} w_i = \frac{b}{\hat{q}}, \quad \hat{q} \partial_{\hat{q}} w_i = b
\]
The second derivative with respect to $\hat{q}$ is:
\[
\partial_{\hat{q}\hat{q}} w_i = -\frac{b}{\hat{q}^2}, \quad \hat{q}^2 \partial_{\hat{q}\hat{q}} w_i = -b
\]
The cross derivative is:
\[
\partial_{q\hat{q}} w_i = 0, \quad q \hat{q} \partial_{q\hat{q}} w_i = 0
\]
From Lemma \ref{lem:optconslingame}\, the maximizer of $\mathcal{H}_{CV}$ is given by:
\[
\alpha^* = \arg \max_{\alpha} \left\{ \ln(\alpha) - \alpha q \partial_q v_i \right\} = (q \partial_q v_i)^{-1}
\]
provided $\partial_q v_i > 0$. Substituting our candidate's derivative $q \partial_q w_i = a$:
\[
\alpha^* = \frac{1}{a}
\]
Substitute $\alpha^* = 1/a$ into $\mathcal{H}_{CV}$:
\[
\begin{aligned}
\mathcal{H}_{CV}(\alpha^*, \partial_q w_i, q) &= \ln\left(\frac{1}{a}\right) + \left(1+\theta _{i}+\frac{\eta _{i}}{N}\right)\ln(q) + \eta_i \ln(\hat{q}) - \left(\frac{1}{a}\right) q \left(\frac{a}{q}\right) \\
&= -\ln(a) + \left(1+\theta _{i}+\frac{\eta _{i}}{N}\right)\ln(q) + \eta_i \ln(\hat{q}) - 1
\end{aligned}
\]
So,
\[
\sup_{\alpha \in \mathbb{R}_+} \mathcal{H}_{CV}(\alpha, \partial_q w_i, q) = -\ln(a) + \left(1+\theta _{i}+\frac{\eta _{i}}{N}\right)\ln(q) + \eta_i \ln(\hat{q}) - 1
\]
We now substitute all the computed parts into the left-hand side (LHS) of the HJB equation \eqref{eqn:HJBgame}, that we repeat here for the reader convenience:
\begin{multline}
\rho v_i\left( q,\widehat{q}\right) -\frac{\widehat{\nu^{-i} }^{2}+\widehat{\xi^{-i}}^{2}}{2}\widehat{q}^{2}\partial_{\widehat{%
q}\widehat{q}}v_i\left( q,\widehat{q}\right) -\frac{\nu _{i}^{2}+\mu
_{i}^{2}}{2}q^{2}\partial_{qq}v_i\left( q,\widehat{q}\right)  -\widehat{q}\partial_{\widehat{q}}v_i\left( q,\widehat{q}\right)
\left( \widehat{g^{-i}}+\frac{\widehat{\nu^{-i} }^{2}+\widehat{\xi^{-i}}^{2}}{2}\right)\\-\gamma q\partial_qv_i(q,\widehat{q})-\sup_{\alpha \in \R_+}\mathcal{H}_{CV}(\alpha,\partial_qv_i,q)- \left(\nu_i \widehat{\nu^{-i}}\right) q \widehat{q} \, \partial_{q\widehat{q}}v_{i}(q, \widehat{q})=0.
\end{multline}%

We have
\[
\rho w_i = \rho (a \ln(q) + b \ln(\hat{q}) + c)
\]
\[
-\frac{\widehat{\nu^{-i} }^{2}+\widehat{\xi^{-i}}^{2}}{2}\widehat{q}^{2}\partial_{\widehat{%
q}\widehat{q}}w_i\left( q,\widehat{q}\right)= -\frac{\widehat{\nu^{-i}}^2 + \widehat{\xi^{-i}}^2}{2} \cdot (-b) = \frac{\widehat{\nu^{-i}}^2 + \widehat{\xi^{-i}}^2}{2} b
\]
\[
-\frac{\nu _{i}^{2}+\mu
_{i}^{2}}{2}q^{2}\partial_{qq}w_i\left( q,\widehat{q}\right)= -\frac{\nu_i^2 + \mu_i^2}{2} \cdot (-a) = \frac{\nu_i^2 + \mu_i^2}{2} a
\]
\[
 -\widehat{q}\partial_{\widehat{q}}w_i\left( q,\widehat{q}\right)
\left( \widehat{g^{-i}}+\frac{\widehat{\nu^{-i} }^{2}+\widehat{\xi^{-i}}^{2}}{2}\right)= -b \left( \widehat{g^{-i}} + \frac{\widehat{\nu^{-i}}^2 + \widehat{\xi^{-i}}^2}{2} \right)
\]
\[
-\gamma q\partial_qw_i(q,\widehat{q})= -\gamma_i \cdot a
\]
\[
-\sup_{\alpha \in \R_+}\mathcal{H}_{CV}(\alpha,\partial_qw_i,q)= - \left( -\ln(a) + \left(1+\theta _{i}+\frac{\eta _{i}}{N}\right)\ln(q) + \eta_i \ln(\hat{q}) - 1 \right) = \ln(a) - \left(1+\theta _{i}+\frac{\eta _{i}}{N}\right)\ln(q) - \eta_i \ln(\hat{q}) + 1
\]
\[
- \left(\nu_i \widehat{\nu^{-i}}\right) q \widehat{q} \, \partial_{q\widehat{q}}w_{i}(q, \widehat{q})= -\nu_i \widehat{\nu^{-i}} \cdot 0 = 0.
\]
Now, summing all these terms together we get the following form for the lefthand side (LHS) of \eqref{eqn:HJBgame}:
\[
\begin{aligned}
\text{LHS} = &\quad \rho a \ln(q) + \rho b \ln(\hat{q}) + \rho c + \frac{\widehat{\nu^{-i}}^2 + \widehat{\xi^{-i}}^2}{2} b + \frac{\nu_i^2 + \mu_i^2}{2} a - b \widehat{g^{-i}} - b \frac{\widehat{\nu^{-i}}^2 + \widehat{\xi^{-i}}^2}{2} \\
&- \gamma_i a + \ln(a) - \left(1+\theta _{i}+\frac{\eta _{i}}{N}\right)\ln(q) - \eta_i \ln(\hat{q}) + 1
\end{aligned}
\]
For $w_i$ to be a solution, the LHS must be zero for all $q, \hat{q} > 0$. Therefore we set coefficient of $\ln(q)$ to zero betting
\[
\rho a - \left(1+\theta _{i}+\frac{\eta _{i}}{N}\right)= 0 \quad \Rightarrow \quad a = \frac{\left(1+\theta _{i}+\frac{\eta _{i}}{N}\right)}{\rho}
\]
and we set coefficient of $\ln(\hat{q})$ to zero getting
\[
\rho b - \eta_i = 0 \quad \Rightarrow \quad b = \frac{\eta_i}{\rho}
\]
These match the definitions given in the statement.

We set the sum of the constant terms to zero getting
\[
\rho c + \frac{\nu_i^2 + \mu_i^2}{2} a - \gamma_i a + \ln(a) + 1 - b \widehat{g^{-i}} = 0
\]

We substitute $a = \left(1+\theta _{i}+\frac{\eta _{i}}{N}\right)/ \rho$ and $b = \eta_i / \rho$
\[
\rho c + \frac{\nu_i^2 + \mu_i^2}{2} \cdot \frac{\left(1+\theta _{i}+\frac{\eta _{i}}{N}\right)}{\rho} - \gamma_i \cdot \frac{\left(1+\theta _{i}+\frac{\eta _{i}}{N}\right)}{\rho} + \ln\left(\frac{\left(1+\theta _{i}+\frac{\eta _{i}}{N}\right)}{\rho}\right) + 1 - \frac{\eta_i}{\rho} \widehat{g^{-i}} = 0
\]

We solve for $\rho c$
\[
\rho c = -\frac{\nu_i^2 + \mu_i^2}{2} \cdot \frac{\left(1+\theta _{i}+\frac{\eta _{i}}{N}\right)}{\rho} + \gamma_i \cdot \frac{\left(1+\theta _{i}+\frac{\eta _{i}}{N}\right)}{\rho} - \ln\left(\frac{\left(1+\theta _{i}+\frac{\eta _{i}}{N}\right)}{\rho}\right) - 1 + \frac{\eta_i}{\rho} \widehat{g^{-i}}
\]

Multiplying both sides by $\rho^{-1}$ and simplify the logarithm term we get

\[
c = -\rho^{-1} \ln(\left(1+\theta _{i}+\frac{\eta _{i}}{N}\right)) + \rho^{-1} \ln(\rho) - \rho^{-1} + \gamma_i \left(1+\theta _{i}+\frac{\eta _{i}}{N}\right)\rho^{-2} - \left(1+\theta _{i}+\frac{\eta _{i}}{N}\right)\rho^{-2} \frac{\nu_i^2 + \mu_i^2}{2} + \rho^{-2} \eta_i \widehat{g^{-i}}
\]

This is exactly the expression for $c$ given in the statement of the present proposition.
    \end{proof}

\subsection*{Proof of Proposition \ref{lem:feedbackgame}}
\begin{proof}
By Lemma \ref{lem:optconslingame} and Proposition \ref{lem:vf1game}
we get that $G_{i}(q,\hat{q})$ is given by
\begin{equation}\label{eqn:c30}
G_{i}(q,\hat{q})=(1+\theta_i+\frac{\eta_i}{N})^{-1}\rho 
\end{equation}
The closed loop equation is the dynamics of  $Q_i$ where we insert \eqref{eqn:c30}, and we find
\begin{eqnarray*}
dQ_i(t)&=&(\gamma_i Q_i(t)-(1+\theta_i+\frac{\eta_i}{N})^{-1}\rho Q_i(t)dt+Q_i(t)\mu_i dW(t)+Q_i(t)\nu_idB(t)\\
Q_i(0)&=&Q_0 \quad \mbox{given}
\end{eqnarray*}
Since the dynamics is linear in $Q_i$,  existence and uniqueness of a solution immediately follows.
\end{proof}

\subsection*{Proof of Theorem \ref{prop:NashN}}
\begin{proof}
First we state the following lemma. The proof works similarly as Proposition $4.3$ of \cite{ML} by taking $N=1,  B=\mu_i q$ and in the simpler case that $\mathcal{L}=\gamma_i$, therefore we do not write the proof here.

\begin{lemma}\label{lem:identity}
      Let $\tau$ be the first exit time from $\R_{++}$.  Then for each $q \in \R_{++}$ and $\alpha \in \mathcal{A}(q)$, we have that
    \begin{equation}\label{eqn:identity}
w_{i}(q)=J(\alpha,q)+\mathbb{E}\left[\int_0^{\tau}e^{-\rho s} [H_{\mbox{\footnotesize{MAX}}}(\partial_qw_{i}(q(s)),q(s))-H_{\mbox{\footnotesize{CV}}}(\alpha(s),\partial_q w_{i}(q(s)),q(s))]ds\right]
    \end{equation}
where
\begin{equation*}
\mathcal{H}_{CV}(\alpha,p,q)=
\ln(\alpha)+\left(1+\theta _{i}+\frac{\eta _{i}}{N}\right)\ln(q)+\eta_i\ln(\widehat{q})-\alpha qp.
\end{equation*}
\[
H_{\mbox{\footnotesize{MAX}}}(p,q)=\sup_{\alpha}\mathcal{H}_{CV}(\alpha,p,q)=-\ln(qp)+\eta_i\ln(\hat{q})-1
\]
\end{lemma}

\noindent We follow the approach of \cite{ML}, Theorem $4.1$.
By Lemma \ref{lem:identity} we have
    $$
    w_{i}(q)\geq J(\alpha,q) \quad \mbox{ for each } \alpha \in \mathcal{A}(q)
$$
since the second term in \eqref{eqn:identity} is positive.
Now notice that, by Proposition \ref{lem:feedbackgame}, we have that  $G^*_{i}(q_0) \in \mathcal{A}(q_0)$ such that
\begin{equation}\label{eqn:Gi2}
G_i(q, \widehat q)=(q\partial_qw_{i})^{-1}=\left(\frac{\left(1+\theta _{i}+\frac{\eta _{i}}{N}\right)}{\rho}\right)^{-1}
=
\frac{\rho}{1+\theta_i+\frac{\eta_i}{N}}.
\end{equation}
Moreover,
by the identity \eqref{eqn:identity} we have
$$
w_{i}(q)=J(G^*_{i}(q_0),q).
$$
Then we deduce \eqref{v_w} and that $G_{i}$ is optimal.
\noindent We finally check that the transversality condition holds, indeed
\begin{eqnarray*}
\lim_{t\rightarrow \infty }\mathbb{E}\left[ e^{-\rho t}v_{i}\left( Q_i(t) \right) \right] &=&%
\lim_{t\rightarrow \infty }e^{-\rho t}\mathbb{E}\left[a \ln(Q_i(t))+b\right]
 \\
&=&a\ln(Q_0)\lim_{t\rightarrow \infty }e^{-\rho t}+a\lim_{t\rightarrow \infty }e^{-\rho t} ( \gamma_i -\alpha ^{\ast }-\frac{_{\nu_i ^{2}+\mu_i ^{2}}}{2})t+\lim_{t\rightarrow \infty }e^{-\rho t}b=0.
\end{eqnarray*}
\end{proof}

\subsection*{Proof of Theorem \ref{thm:max_principle}}
\begin{proof}
     The necessary part of the proof is a straightforward adaptation of the Maximum principle presented in \cite[Theorem 5.19]{carmona2016lectures}.\\
     The sufficient part is a consequence of Arrow condition, indeed
     \begin{align*}
H^i_{Max}(Q,Y^i,Z_i)&=\sup_\alpha H^i(Q,Y^i,Z_{i},\alpha)\\
&=\sum_{k=1}^N (\gamma_kQ_k(t)-\frac{1}{Y^{k,k}(t)})\cdot Y^{i,k}(t)-\log(Q_iY^{i,i}(t))+(1+\theta_i)\log(Q_i)+\frac{\eta_i}{N}\sum_{j=1}^N\log(Q_j)+\\
&\hspace{6cm}+\sum_{k=1}^N Z_{i,k,k}(t)\mu_k Q_k+\sum_{k=1}^N Z_{i,k,0}(t)\nu_kQ_k\\
\end{align*}
and therefore
\[\partial^2_{Q_iQ_i}H^i_{Max}=-\left(\theta_i+\frac{\eta_i}{N}\right)\frac{1}{Q_i^2}, \quad \partial^2_{Q_jQ_j}H^i_{Max}=-\left(\frac{\eta_i}{N}\right)\frac{1}{Q_j^2},\quad \partial^2_{Q_iQ_j}H^i_{Max}=-0,\]
    which imply that $H^i_{Max}(Q,Y^i,Z_i)$ is jointly concave in the state variables.
\end{proof}

\subsection*{Proof of Proposition \ref{prop:open_loop}}
  \begin{proof}
 In order to find open loop equilibria we will proceed in the following way.
\begin{itemize}
    \item Formulate an Ansatz
    \item Plug in the equation for $Y^{i,j}$ the ansatz.
    \item Derive an equation for $Y^{i,j}$ from the ansatz, by Ito Formula
    \item Compare the two results to characterize the ansatz and conclude by deriving the open loop equilibria through $Y^i$.
\end{itemize}
The ansatz is
\[Y^{i,j}(t)=\phi^i\frac{1}{Q^i(t)}\delta_{ij}+\psi^{j}\frac{1}{Q^j(t)}(1-\delta_{ij})\]

By plugging the ansatz in the equation for $Y^{i,j}$ we get

\begin{multline}\label{eq:Y_1}
    dY^{i,j}(t)=\rho \left(\phi^i\frac{1}{Q^i(t)}\delta_{ij}+\psi^{j}\frac{1}{Q^j(t)}(1-\delta_{ij})\right)-\left((\gamma_j-\alpha_j)\phi^i+(1+\theta_i)\right)\frac{1}{Q^i(t)}\delta_{ij}dt+\\
    \hspace{5cm}-\left((\gamma_j-\alpha_j)\psi^j+\frac{1}{N}\eta_j\right)\frac{1}{Q^j(t)}(1-\delta_{ij})dt-Z_{i,j,j}(t)\mu_jdt-Z_{i,j,0}(t)\nu_jdt\\
    \hspace{10cm}+\sum_{k=0}^N Z_{i,j,k}(t) dW^k(t)
\end{multline}
By applying Ito Formula on the ansatz,

\begin{multline}\label{eq:Y_2}
   dY^{i,j}(t)=-\left(\phi^i_t(\gamma_i-\alpha_i)-(\phi^i_t)'+\mu_i^2+\nu_i^2\right)\frac{\delta_{ij}}{Q^i(t)}dt-\left(\psi^j_t(\gamma_j-\alpha_j)-(\psi^j_t)'+\mu_j^2+\nu_j^2\right)\frac{1-\delta_{ij}}{Q^j(t)}dt+\\
   \hspace{2cm}+\mu_i\frac{\delta_{ij}}{Q^i(t)}dW^i(t)+\mu_j\frac{1-\delta_{ij}}{Q^j(t)}dW^j(t)+\left(\nu_j\frac{\delta_{ij}}{Q^j(t)}+\nu_i\frac{1-\delta_{ij}}{Q^i(t)}\right)dW^0(t)\\
\end{multline}

By comparing \eqref{eq:Y_1} and \eqref{eq:Y_2}, we get
\[Z_{i,j,0}=\left(\nu_j\frac{\delta_{ij}}{Q^j(t)}+\nu_i\frac{1-\delta_{ij}}{Q^i(t)}\right), \quad Z_{i,j,j}=\left(\mu_j\frac{1-\delta_{ij}}{Q^j(t)}\right)\]
\[Z_{i,j,i}=\left(\mu_i\frac{\delta_{ij}}{Q^i(t)}\right)\]
and
\[(\phi^i)'=\rho\phi^i-(1+\theta_i),\quad (\psi^j)'=\rho\psi^j-\frac1N \eta_j.\]
The transversality condition gives the following final condition on the coefficient,
\[\lim_{t\to\infty}e^{-\rho t}Y^{i,j}_tQ^i_t=\lim_{t\to\infty}e^{-\rho t}\left(\phi^i(t)\delta_{ij}+\psi^j(t)(1-\delta_{ij})\right)=0\]
Therefore, the solution of the system
\[\begin{cases}
(\phi^i)'=\rho\phi^i-(1+\theta_i)\\
\lim_{t\to\infty}e^{-\rho t}\phi^i(t)=0
\end{cases}\]
is
\[\phi^i=\frac{1+\theta_i}{\rho}.\]
And similarly $\psi^j=\frac{\eta_j}{N\rho}$.
From the optimality condition, we conclude
\[\alpha^i_{OL}(t)=\frac{1}{Q^i(t)Y^{i,i}(t)}=\frac{1}{\phi^i(t)}=\frac{\rho}{1+\theta_i}\]
    \end{proof}

\subsection*{Proof of Proposition \ref{prop:Price_of_anarchy}}
\begin{proof}
    Observe that if $\theta_i\equiv\theta$ and $\eta_i\equiv\eta$, we have
    \[\alpha^i_{SP}\equiv\alpha_{SP}, \quad \alpha^i_{CL}\equiv\alpha_{CL}.\]
    Therefore, by evaluating $J_{SP}$ in $\boldsymbol{\alpha}_{SP}$ defined in \eqref{alpha_SP} and $J^i$ in $\alpha_{CL}$ defined in \eqref{alpha_SP}, we have
    \[J_{SP}\left(\boldsymbol{\alpha}_{SP}\right)=\frac{1}
    {\rho}\log(\alpha_{SP})+\frac{1+\theta+\eta}{\rho^2}\left(\gamma-\alpha_{SP}-\frac12\left(\mu^2+\nu^2\right)\right),\]
    \[J^i\left(\alpha_{CL}\right)=\frac{1}{\rho}\log(\alpha_{CL})+\frac{1+\theta+\eta}{\rho^2}\left(\gamma-\alpha_{CL}-\frac12\left(\mu^2+\nu^2\right)\right).\]
    Since the function $F(x)=\frac{1}{\rho}\log(x)+\frac{1+\theta+\eta}{\rho}\left(\gamma-x-\frac12\left(\mu^2+\nu^2\right)\right)$ is decreasing for $x\geq \alpha_{SP}$ and since $\alpha_{SP}<\alpha_{CL}$, we conclude that $J_{SP}(\boldsymbol{\alpha}_{SP})-J^i(\alpha_{CL})>0$.
\end{proof}

\subsection*{Proof of Theorem \ref{thm:convme}}
\begin{proof}
    We start exploiting the definition of $w^{N,i}$. From \eqref{vni_explicit} we get
\begin{align}\label{wni}
w^{N,i}(q_i,m)&=\int\limits_{\R^{N-1}}w_i( q, \hat q)\prod\limits_{j\neq i}m(dq_j)\\
&=\int\limits_{\R^{N-1}}\left( a\ln(q)+\tilde c\right)\prod\limits_{j\neq i}m(dq_j)
+\int\limits_{\R^{N-1}}b\ln(\widehat q)\prod\limits_{j\neq i}m(dq_j)\,.
\end{align}
Observe that
$$
\int\limits_{\R^{N-1}}\left( a\ln(q)+\tilde c\right)\prod\limits_{j\neq i}m(dq_j)=a\ln(q)+\tilde c\,,
$$
since $m$ is a probability measure and the integrand does not depend on $q_j$ for $j\neq i$.
Analyzing the second term in \eqref{wni}, we get
\begin{align*}
\int_{\R^{N-1}}\ln(\hat q)\prod\limits_{j\neq i}m(dq_j)&=\frac 1{N}\sum\limits_{k\neq i}\int_{\R^{N-1}}\ln(q_k)\prod\limits_{j\neq i}m(dq_j)\\
&=\frac 1{N}\sum\limits_{k\neq i}\int_{\R^{N-2}}\left(\int_\R \ln(q_k)\, m(dq_k)\right)\prod_{j\neq i,k}m(dq_j)\\
&=\frac 1{N}\sum\limits_{k\neq i}\int_\R \ln(q_k)\,m(dq_k)=\frac{N-1}{N}\int_{\R_+} \ln(q)\,m(dq)\,.
\end{align*}
This means that we have
$$
w^{N,i}(q_i,m)=a\ln(q)+\tilde c+\frac{b(N-1)}{N}\int_{\R_+} \ln(q)\,m(dq)
$$
Then we can write
\begin{eqnarray*}
\int_{\R_+} \big|w^{N,i}(q,m)-U(q,m)\big|\,m(dq)
&\leq & \int_{\R_+}\left|\tilde c-c+b\left(\frac{(N-1)}{N}\int_{\R_+} \ln(q)\,m(dq)-\int_{\R_+}\ln(q) m(dq)\right) \right|\, m(dq)\\ &=&\int_{\R_+}\left|\tilde c-c+\frac{b}{N} \int_{\R_+}\ln(q) m(dq)\right|\, m(dq)\\ &=& \frac{1}{N}\left(\rho^{-2}\eta\left(\gamma-\frac{\rho}{1+\theta}-\frac{\mu^2}{2}\right)+b \left| \int_{\R_+}\ln(q) m(dq)\right|\right)
\end{eqnarray*}
\end{proof}
\newpage

\bibliographystyle{plain}
\bibliography{biblio}

\end{document}